 \documentclass[11pt]{article}
\usepackage{hyperref}
\usepackage{cite,amsfonts,amssymb,amsmath,exscale,bbm}
\usepackage{footnpag} 
\usepackage[all]{xy}  

 \usepackage{bm}

\usepackage[mathscr]{eucal}
\usepackage{amsfonts,amssymb,amsmath,bbm}
\usepackage{bbm}

\usepackage[all]{xy}  
\usepackage{tikz}

\newcommand{\beq}{\begin{equation}}
\newcommand{\eeq}{\end{equation}}
\newcommand{\be}{\begin{equation}}
\newcommand{\ee}{\end{equation}}
\newcommand{\beqa}{\begin{eqnarray}}
\newcommand{\eeqa}{\end{eqnarray}}
\newcommand{\bean}{\begin{eqnarray*}}
\newcommand{\eean}{\end{eqnarray*}}

\newcommand{\cP}{{\mathcal P}}

\newcommand\U{{\mathrm U}}

\renewcommand\P{{\mathcal{P}}} 

\newcommand{\R}{\mathbb{R}}
\newcommand{\C}{\mathbb{C}}

\newcommand{\Z}{\mathbb{Z}}

\newcommand{\SO}{\mathrm{SO}}
\newcommand{\SU}{\mathrm{SU}}

\def\maps{\colon}

\def\to{\rightarrow}

\def\direct{\int^{\oplus}}
\def\extd{\mathrm {d}}

\def\s{\section}
\def\ss{\subsection}
\def\sss{\subsubsection}

\newcommand\unit{\mathbbm{1}}

\newtheorem{lemma}{Lemma}

\newcommand{\urtriangle}{ \begin{tikzpicture}
\draw (0,0.15) -- ++(0.15,0) -- ++(0,-0.15) --cycle ;
\end{tikzpicture}}

\newcommand{\ultriangle}{\begin{tikzpicture}
\draw (0,0) -- ++(0,0.15) -- ++(0.15,0) --cycle ;
\end{tikzpicture}}

\newcommand{\drtriangle}{\!   \begin{tikzpicture}
\draw (0,0) -- ++(0.15,0) -- ++(0,0.15) --cycle ;
\end{tikzpicture}}

\newcommand{\dltriangle}{\begin{tikzpicture}
\draw (0,0) -- ++(0,0.15) -- ++(0.15,-0.15) --cycle ;
\end{tikzpicture}}

\newcommand{\tinyurtriangle}{ \begin{tikzpicture}[scale=0.7]
\draw (0,0.15) -- ++(0.15,0) -- ++(0,-0.15) --cycle ;
\end{tikzpicture}}

\newcommand{\tinyultriangle}{\begin{tikzpicture}[scale=0.7]
\draw (0,0) -- ++(0,0.15) -- ++(0.15,0) --cycle ;
\end{tikzpicture}}

\newcommand{\tinydrtriangle}{\!   \begin{tikzpicture}[scale=0.7]
\draw (0,0) -- ++(0.15,0) -- ++(0,0.15) --cycle ;
\end{tikzpicture}}

\newcommand{\tinydltriangle}{\begin{tikzpicture}[scale=0.7]
\draw (0,0) -- ++(0,0.15) -- ++(0.15,-0.15) --cycle ;
\end{tikzpicture}}

\newcommand{\downtriangle}{ \begin{tikzpicture}]
\draw (0,0) -- ++(0.15,0) -- ++(-0.075,-0.15) --cycle ;
\end{tikzpicture}}

\newcommand{\squarediagup}{\begin{tikzpicture}
\draw (0,0) rectangle (0.15, 0.15);
\draw (0,0) -- (0.15,0.15);
\end{tikzpicture}}

\newcommand{\squarediagdown}{\begin{tikzpicture}
\draw (0,0) rectangle (0.15, 0.15);
\draw (0.15,0) -- (0,0.15);
\end{tikzpicture}}

\newcommand{\squarediag}{\begin{tikzpicture}
\draw (0,0) rectangle (0.15, 0.15);
\draw (0.15,0) -- (0,0.15);
\draw (0,0) -- (0.15,0.15);
\end{tikzpicture}}

\newcommand{\squarenodiag}{\begin{tikzpicture}
\draw (0,0) rectangle (0.15, 0.15);
\end{tikzpicture}}

\newcommand{\pentagon}{\begin{tikzpicture}[scale=0.1]
\draw (-0.5, 0) --(-1, 1.1) -- (0, 1.7) -- (1,1.1) -- (0.5, 0) -- (-0.5, 0);
\end{tikzpicture}}

\begin{document}




\centerline{\Large \bf A 2-categorical  state sum model}

\bigskip

\centerline{Aristide Baratin$^{1}$, Laurent Freidel$^{2}$}

\bigskip 

{\small
\centerline{${}^1$ Department of Applied Mathematics, University of Waterloo } 
\centerline{200 University Ave W, Waterloo ON, N2L 3G1, Canada}
\vskip .3em
\centerline{${}^2$ Perimeter Institute for Theoretical Physics} 
\centerline{31 Caroline Str. N, 
Waterloo ON, N2L 2Y5, Canada}

\medskip

\begin{abstract}  
It has long been argued that higher categories provide the proper algebraic structure underlying state sum invariants of 4-manifolds.
This idea has been refined recently, by proposing to use 2-groups and their representations as specific examples of 2-categories. The challenge has been to make these proposals fully explicit. Here we give a concrete realization of this program. Building upon our earlier work with Baez and Wise on the representation theory of 2-groups, we construct a four-dimensional state sum model based on a categorified version of the Euclidean group. We define and explicitly compute the simplex weights, which may be viewed a categorified analogue of Racah-Wigner 6$j$-symbols. These weights solve an hexagon equation that encodes the formal invariance of the state sum under the Pachner moves of the triangulation.  This result unravels the combinatorial formulation  of the Feynman amplitudes of quantum field theory  on flat spacetime  proposed in \cite{BaratinFreidel}, which was shown to lead after gauge-fixing to Korepanov's invariant of 4-manifolds.
\end{abstract}

\medskip

\s{Introduction}

This paper results from the  convergence of two lines of investigation in state sum models: 
state sums can be viewed either as  topological objects \cite{BaratinFreidel}  or categorical constructions \cite{BarrettMackaay}.
Well understood in dimension two \cite{FHK} and three \cite{PonzanoRegge,TuraevViro}, this convergence is established here in dimension four.

State sum models provide a powerful technical tool for the combinatorial construction of manifold invariants and topological quantum field theories. The idea is to rewrite a path integral as a sum of local weights defined using a triangulation of the manifold, and to reformulate topological invariance as a set of algebraic equations for the weights. These equations encode the invariance of the state sum under elementary re-buildings of the triangulation, the so-called Pachner moves,  which are known to relate any topologically equivalent configurations. 
This procedure is the core of the lattice definition of two-dimensional topological field theory by Fukuma, Hosono and Kawai \cite{FHK}. 
Notorious examples of state sum models in three dimensions are the Ponzano-Regge and Turaev-Viro models based on the representation category of a (quantum) group, leading to a state sum formulation of quantum gravity in three-dimensional spacetime  \cite{PonzanoRegge, TuraevViro}.  

State sum models are also at the root of the spin foam approach to quantum gravity in four dimensions \cite{Perez}. Stemming from a formulation of gravity as a constrained topological theory, the main strategy in this approach has been to quantize the topological theory using a state sum model and to impose the constraints in the resulting quantum theory. Specific realizations of this idea led to a background independent formulation of the gravity path integral as a sum over geometries  displaying a fundamental discreteness at the Planck scale.  Remarkably, this formulation enables one to define transition amplitudes between the states of the gravitational field in loop quantum gravity \cite{FK, EPRL}. 
The close relation to topological field theory is one of the striking features of this approach: not only is it sufficient to determine the form of the boundary states, but it may arguably enable one to keep the diffeomorphism symmetry, and with it the low energy behavior of the theory, under control. 

Category theory appears to be a natural arena for constructing state sum models and generalizing the Fukuma, Hosono and Kawai's procedure to  dimensions higher than two \cite{BaezDolan, 4DTQFT, BaezLauda}. This is due to a remarkable correspondence between the combinatorics of Pachner moves and the coherence laws of (higher) categories. Well understood in three dimensions \cite{BarrettWestbury}, this correspondence is however much more difficult to exploit in higher dimensions, due to the complexity of higher algebraic structures.
Formalisms have been proposed for the construction of four-dimensional models using 2-categories \cite{Mackaay}, 
but the only known examples of such models use a very restricted class of 2-categories, ones with a single object \cite{CKY, Mackaay2}. 

As an attempt to find analogues of three-dimensional models built from representations of a group, Barrett and Mackaay proposed in \cite{BarrettMackaay} to build four-dimensional state sum models starting with a categorical group, or 2-group. 2-groups play the same role in higher gauge theory as groups do in gauge theory: just as groups can be used to describe connections defining parallel transport along curves, 2-groups can be used to describe `2-connections'  defining parallel transport along both curves and surfaces \cite{BaezHuerta}. 
 A 2-group can be defined as a `crossed module', which is a pair of groups related by an homomorphism $H \to G$ and an action of $G$ on $H$ satisfying some compatibility conditions. One of the simplest examples of (Lie) 2-group is the Poincar\'e 2-group \cite{Baez}, determined  by the Lorentz group $G$ acting on the translation group $H$ of Minkowski space, taken with the trivial homomorphism $H \to G$. The idea of \cite{BarrettMackaay}, also expressed and further explored by Crane, Sheppeard and Yetter in \cite{CraneYetter2, CraneYetter, CraneSheppeard},  was to try using the representations of the Poincar\'e 2-group to construct a four-dimensional analogue of the Ponzano-Regge model for three-dimensional quantum gravity.  

This is the first of the two lines of investigation leading to the results presented here. We explicitly construct the model outlined in these works, using a 2-group closely related to the Poincar\'e 2-group, the Euclidean 2-group. The core of our construction is the complete calculation of the weight for the 4-simplex, which can be seen as a categorified analogue of Racah-Wigner 6$j$-symbols. The Euclidean 2-group is a categorified version of the Euclidean group: it differentiates the roles of rotations and translations by treating the former as objects in a category and the latter as morphisms. Consequently its representation theory looks quite different from that of the group \cite{2Rep}. In fact, just as the representations of a group can be viewed as objects in a category, the representations of a 2-group can be viewed as objects in a 2-category. The model developed in this paper thus gives an explicit realization of the 2-categorical approach to state sum models in four dimensions. 

The second line of investigation concerns the combinatorial reformulation of Feynman amplitudes in quantum field theory on flat spacetime proposed in \cite{BaratinFreidel3d, BaratinFreidel}. 
The goal of that work, motivated by earlier results in three-dimensional spin foam gravity \cite{Barrett, PR3a, PR3b},  was to try to bridge the gap between the algebraic framework of spin foams  and the standard formalism of quantum field theory formulated on a background spacetime with a fixed metric.  The result was to show that Feynman amplitudes on flat space-time can be rewritten as Wilson line observables in a state sum model characterized by specific weights. 
In this approach the state sum replaces entirely the background geometry,  which arises only after the state sum is performed.
Moreover, when evaluated on closed manifolds without Wilson line insertion, and after gauge-fixing,  the state sum model was shown in \cite{BaratinFreidel}  to reproduce the 4-manifold invariant previously constructed by Korepanov \cite{Korepanov1, Korepanov2, Korepanov3}.

The original motivation for the present paper was to unravel the algebraic nature of the state sum model discovered in \cite{BaratinFreidel}. 
This is where the two lines of investigation meet: as we show here,  the relevant structure is the 2-category of representations of the Euclidean 2-group. 
The fact that such a 2-categorical model shows up naturally  in ordinary quantum field theory is an intriguing and exciting outcome of our work. 

Explicitly, given a  triangulated closed 4-manifold $\Delta$, the state sum model of \cite{BaratinFreidel} is characterized by weights which are (real) functions 
of a set of positive numbers  $l_e \!\in\! \R_+$ labeling the edges $e$ and a set of integer spins $s_t \!\in\! \Z$ labeling the triangles $t$. These weights are given by the formula:
 \beq \label{weights}
W_{\Delta} (l_e, s_t)  = \prod_{t \in \Delta} {2 {A}_t(l_e)} \prod_{\sigma \in \Delta}\, \frac{\cos \left[\sum_{t\subset \sigma}\! s_t \phi^\sigma_t(l_e)\right]}{V_\sigma(l_e)} 
\eeq
where the products are over all triangles $t$ and all 4-simplices $\sigma$. 
There is a factor $2 {A}_t$ for each triangle $t$, which depends on the three labels $l_e \!\in\! \R_+$ on the edges of $t$.  Whenever these numbers satisfy the triangle inequality, i.e if they define a Euclidean geometry for $t$, $A_t$ is equal to the area of the triangle; otherwise,  it is equal to zero.  
This means the weight $W_{\Delta}$ is zero unless the set of labels $l_e$ is consistent with all triangle inequalities. There is a factor for each 4-simplex $\sigma$, which depends on ten labels $l_e\!\in\! \R_+$,  one for each edge of $\sigma$; and ten labels $s_t\!\in \!\Z$, one 
 for each of  triangle of $\sigma$. $V_\sigma(l_e)$ is  equal to $4!$ times the volume of the 4-simplex with edge lengths $l_e$; the sum in the argument of the cosine is over all ten triangles of $\sigma$, and $\phi_t^\sigma$  denotes the dihedral angle between the two tetrahedra of $\sigma$ that meet at the triangle $t$.  The partition function, or state sum, is formally obtained by summing the weights (\ref{weights}) over all values of the labels, using the Lebesgue measures $\extd l_e^2$ for the real variables and the counting measure $\frac{1}{2\pi} \!\sum_{s_t}$ for the integer variables. 

The conjectured 2-categorical nature of the weights (\ref{weights})  prompted the in-depth study \cite{2Rep} of the infinite dimensional representation theory of Lie 2-groups introduced in \cite{CraneYetter}. The case of the Euclidean 2-group, determined by the pair of groups $G\!=\!\SO(4)$ and $H\!=\!\R^4$, was described in \cite{Baratin-Wise}. In Section \ref{labelling} we begin by reviewing this description. We explain how the variables $l_e \!\in\! \R_+$ and $s_t \!\in\! \Z$ can be understood as labeling  `irreducible representations' of the Euclidean 2-group and  `irreducible 1-intertwiners' between representations. A key aspect of this theory is that any such 1-intertwiner defines an ordinary representation of the rotation group SO(4). There are also `2-intertwiners' between 1-intertwiners, which define ordinary intertwiners between SO(4) representations.  So, just like models built from SO(4) representation theory, the model based on the representation theory of the Euclidean 2-group starts with an assignment of a representation of SO(4) on each triangle and an SO(4) intertwiner on each tetrahedron. The crucial difference, however, is that in the latter case the representation on each triangle is infinite-dimensional and depends on the data labeling the bounding edges. 

Section \ref{symbol} presents our main result: the definition and computation of the weights of the state sum. 
The weight for each 4-simplex, coined `10$j$ 2-symbol' in this paper, gives an analogue of Racah-Wigner 6$j$ symbols in the 2-categorical context.  
We explain the construction in explicit detail and perform a full calculation of the weights, which leads precisely to the formula (\ref{weights}).  

We conclude with an outlook in Section \ref{outlook}.

\s{Labeled triangulations}  \label{labelling} 

In this section we review the key ingedients from the representation theory of the Euclidean 2-group required for the construction of the state sum model. 

Given an triangulated 4-manifold $\Delta$, we will: 
\begin{itemize}
\item assign an irreducible representation of the Euclidean 2-group to each edge of $\Delta$
\item assign an irreducible 1-intertwiner to each triangle of $\Delta$, and
\item assign a 2-intertwiner to each tetrahedron.
\end{itemize}
The state sum model gives a way to compute an amplitude for any such assignment. 
We begin with a geometrical description of what the notions listed above amount to.  We refer the reader to \cite{2Rep} for the full details of the representation theory.

\ss{Representations on edges} 

A representation of the Euclidean 2-group is given by an  $\SO(4)$-equivariant map $\chi\maps X \to \R^4$, where $X$ is some space on which the rotation group acts, $(g, x) \mapsto g\cdot x$.   `Equivariant' means that $\chi$ commutes with the group action:
\beq
\chi(g\cdot x) = g\,\chi(x)
\eeq
where $g\chi(x)$ is the image of the $\R^4$-vector $\chi(x)$ by the rotation $g\in \SO(4)$.  Irreducible representations are the ones for which the action on $X$ is transitive and the map $\chi$ is one-to-one: in this case  $X$ is identified to a single $\SO(4)$-orbit in $\R^4$, a 3-sphere of given radius $l \in \R_+$:
\beq
X_l = \{x\in \R^4, |x|=l\}
\eeq
So, irreducible representations are effectively labeled by positive numbers, the radii of spheres. There is a notion of tensor product of representations \cite{2Rep}: the tensor product $l\otimes m$ of two irreducible representations corresponds to the map $X_l \times X_m \to \R^4$ given by 
$(x,x')\mapsto x+ x'$.

In what follows, we will work with the spherical measures $\extd^3_{l} x$ induced on $X_l$ by the Lebesgue measure on $\R^4$ and normalized as:
\beq \label{MeasSp}
\extd^3_{l} x  := \frac{1}{\pi} \extd^4 x \, \delta(|x|^2 - l^2).
\eeq
giving a total volume of $\pi l^3$ for the 3-sphere $X_l$. 

Other notations are as follows. We assume $\R^4$ is equipped with its standard 
basis of unit vectors  $e_1, e_2, e_3$ and $e_4$. We identify U(1) with the subgroup of SO(4) rotations that leave the whole plane $(e_1, e_2)$ fixed; we denote by $h_{\phi}$ the U(1) element rotating the plane $(e_3, e_4)$ by the angle $\phi$.  In Section \ref{2int} below we will consider the SO(3)-subgroup of rotations that leave $e_4$ fixed; abusing notation, we will refer to this subgroup simply as SO(3).

\ss{1-Intertwiners on triangles} 
\label{1-int}

Consider three irreducible representations labeled by $l,m,n \in \R_+$ satisfying the  (strict) triangle inequality. 
A 1-intertwiner between $l\otimes m$ and $n$, drawn as 
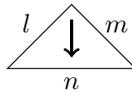
\[
\begin{array}{c}
\begin{tikzpicture}
\draw (0,0) -- (0.85,0.85) -- (1.7, 0) --cycle;
\path (0.25,0.6) node {\small $l$};
\path (1.45,0.55) node {\small $m$};
\path (0.85,-0.2) node {\small $n$};
\draw[->,very thick]
(0.85, 0.65) -- (0.85, 0.15);
\path (0.70, 0.37); 
\end{tikzpicture}
\end{array}
\]
is described as follows. Let $T$ be the set of triples of $\R^4$-vectors $(x, y, z)$  forming a triangle of lengths $l,m,n$:
\beq \label{triangle}
|x| = l, \, |y|=m, \, |z|=n, \quad x+y=z
\eeq
It will be convenient to use the symbolic notation `$\vartriangle$'  for a generic point $(x, y, z)$ in $T$. 
The natural action of SO(4) on $\R^4$ induces a transitive action on $T$, which we write $(g, \vartriangle) \mapsto g \!\!\vartriangle$. 
A 1-intertwiner amounts to a SO(4) vector bundle over $T$, that is, a vector bundle $V$ with a fiber-preserving $\SO(4)$ action. 
More precisely, for $\vartriangle\in T$, $\varphi \in V_{\!\vartriangle}$, and $g \in \SO(4)$, we can write the action as 
\beq \label{SO(4)rep}
g (\!\vartriangle, \varphi)  = (g\!\!\vartriangle, \Phi^g_{\!\vartriangle}(\varphi)),
\eeq where $\Phi^g_{\!\vartriangle}\maps V_{\!\vartriangle}\to V_{\! g \!\vartriangle}$ are invertible linear maps satisfying the rule:
\begin{equation}
\label{cocycle}
\Phi^{gg'}_{\!\vartriangle} = \Phi^g_{g' \!\vartriangle} \Phi^{g'}_{\!\vartriangle}
\end{equation}
This rule says that (\ref{SO(4)rep}) defines a representation of $\SO(4)$ on {\em sections} of the vector bundle $V$.  Alternatively, if we fix a  point $\!\vartriangle \in T$ and restrict the action to its stabilizer 
\beq
G_{\!\vartriangle} := \{h \in \SO(4): h\!\!\vartriangle  = \vartriangle\}
\eeq
then (\ref{cocycle}) is just the equation for a representation of $G_{\!\vartriangle}$ on $V_{\!\vartriangle}$. 
As shown in \cite{2Rep}, two such 1-intertwiners defining equivalent $G_{\!\vartriangle}$ representations at the same point $\!\vartriangle \in T$ belong to the same equivalence class
In the language of Mackey  \cite{Mackey}, the full $\SO(4)$ representation is the one {\it induced} by the $G_{\!\vartriangle}$ representation.  Irreducible 1-intertwiners are the ones for which the $G_{\!\vartriangle}$ representation is irreducible. Here the SO(4) action on $T$ has a U(1) subgroup as stabilizer group: its irreducible representations are one-dimensional, so that the bundle $V$ is a line bundle $V_{\vartriangle} \simeq \C$. Irreducible 1-intertwiners are labeled by elements of Irrep(U(1)) $\simeq \Z$. 

Since SO(4) acts transitively on $T$, the rule (\ref{cocycle})  determines all maps $\Phi^{g}_{\!\vartriangle}$ in terms of the map at a given point $\vartriangle \in T$. For practical calculations, which will require to pick a representative in each equivalence class of 1-intertwiners, it will be convenient  to single out a `reference'  triangle in $T$. For each set of values for the edge labels, 
we denote by $\vartriangle^{\!\! o}$ the unique triangle $(x_0, y_0, z_0) \in T$ such that:
\beq  \label{reftriangle}
x_0= l e_1, \quad 
z_0=n(\sin \gamma e_1+\cos \gamma e_2)
\eeq 
for some $\gamma \in [0,\pi]$, where $e_i$ are the basis vectors. Thus, the pair of vectors $(x_0, z_0)$  defines a triangle with positive orientation in the oriented plane $(e_1, e_2)$, with U(1) as stabilizer group $G_{\!\vartriangle^{\!o}}$. The maps $\Phi^{g}_{\!\vartriangle}$ at any point $\vartriangle = k\!\! \vartriangle^{\!\!o}$, where $k\in \SO(4)$,  can then be expressed as:
\beq \label{rule}
\Phi^{g }_{k\! \vartriangle^{\!o}}  =\Phi^{gk}_{\vartriangle^{\!o}}  \left(\Phi^k_{\vartriangle^{\!o}}\right)^{-1}  
\eeq
For each $s\in \Z$, let us fix once and for all a nowhere vanishing complex function $g \mapsto \Phi^{s}(g)$ on SO(4) such that $\Phi^{s}(1) = 1$ and
%
\beq \label{spincocycle}
\Phi^{s}(g h_{\theta})=  e^{is\theta} \Phi^{s}(g), 
\eeq
for all $g\!\in\! \SO(4)$ and all U(1) rotation $h_\theta$ of angle $\theta$. 
A representative in the equivalence class of irreducible 1-intertwiners labeled by $s$ is obtained by choosing $\Phi^{g }_{\vartriangle^{\!o}}\maps V_{\!\vartriangle^{\!o}}\to V_{\! g \!\vartriangle^{\!o}}$, which is a map $\C \to \C$,  to act by multiplication by the complex number $\Phi^{s}(g)$. Note that the set of such functions (\ref{spincocycle}) is not empty: consider for example the product
\[
\Phi^{s}(g) := D^{j}_{\frac{s}{2}\frac{s}{2}}(g_L) \, D^{j}_{\frac{s}{2}\frac{s}{2}}(g_R)
\] of two $\SU(2)$ Wigner $D$-matrices diagonal elements in a given  irreducible SU(2) representation  $j$,  where $g_L,g_R$ are the left and right SU(2) components of $g$  and U(1) is identified to the subgroup of diagonal elements $g_L \!=\! g_R \!:=\! h_{\theta}$.  

%

Using an $\SO(4)$-invariant measure $\mu_T$ on $T$, we may also view the SO(4) representation (\ref{SO(4)rep}) as acting on the space
\beq \label{directint}
\direct \extd\mu_{T}(\vartriangle) V_{\!\vartriangle}
\eeq
of $L^2$ sections of the bundle $V$  -- provided  the vector bundle is also a  `Hilbert bundle', namely the fibers $V_{\vartriangle}$ are Hilbert spaces. When the 1-intertwiner is irreducible,  $V$ is a line bundle; in this case the direct integral is just the space $L^2(T, \mu_{T})$ of complex functions $f$ for which
\beq 
\int_T \extd \mu_T(\vartriangle) |f(\vartriangle)|^2 < \infty.
\eeq 
Since SO(4) acts transitively on $T$, all invariant measures coincide up to a multiplicative constant.  
We will use the following measure on $T$:
\beq \label{MeasInt}
\extd\mu_{T} =   \extd^3_l x \, \extd^3_m y \, \extd^3_n z \, \delta^4(x+y-z)
\eeq
expressed in terms of the spherical measures (\ref{MeasSp}). With the chosen normalization, the total volume of $T$ gives the  area $A$ of a triangle
of lengths $l,m,n$:
\beq \label{area}
\mu_T(T) = 2  A(l,m,n).
\eeq 


The composition of two 1-intertwiners 
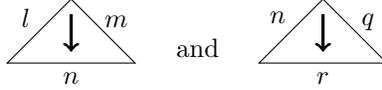
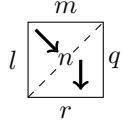
\[
\begin{array}{c}
\begin{tikzpicture}
\draw (0,0) -- (0.85,0.85) -- (1.7, 0) --cycle;
\path (0.25,0.6) node {\small $l$};
\path (1.45,0.55) node {\small $m$};
\path (0.85,-0.2) node {\small $n$};
\draw[->,very thick]
(0.85, 0.65) -- (0.85, 0.15);
\path (0.70, 0.37); 
\end{tikzpicture}
\end{array}
\quad \mbox{and} \quad 
\begin{array}{c}
\begin{tikzpicture}
\draw (0,0) -- (0.85,0.85) -- (1.7, 0) --cycle;
\path (0.25,0.6) node {\small $n$};
\path (1.45,0.55) node {\small $q$};
\path (0.85,-0.2) node {\small $r$};
\draw[->,very thick]
(0.85, 0.65) -- (0.85, 0.15);
\path (0.70, 0.37); 
\end{tikzpicture}
\end{array}
\]
drawn as: 
\[
\begin{array}{c}
\begin{tikzpicture}
\draw (0,0) rectangle (1, 1);
\path (-0.2, 0.5) node {\small $l$};
\path (0.5, 1.2) node {\small $m$};
\path (1.15, 0.5) node {\small $q$};
\path (0.5, -0.2) node {\small $r$};
\path (0.42, 0.83); 
\path (0.88, 0.37); 
\path (0.5, 0.5) node {\small $n$};
\draw[dashed] (0,0) -- (0.45,0.45);
\draw[dashed] (0.6, 0.6) -- (1, 1);
\draw[->, very thick]
(0.1,0.9) -- (0.4, 0.6);
\draw[->, very thick]
(0.7, 0.5) -- (0.7, 0.1);
\end{tikzpicture}
\end{array}
\]
is described as follows. Consider the set $Q$ of quadruples of $\R^4$-vectors $(x, y, v, w)$ forming a quadrangle of lengths $l, m, q, r$:
\beq \label{quadra}
|x| = l, \, |y|=m, \, |v| = q, \, |w|=r,\, \quad x+y+v =w
\eeq
Let $Q^-_n \subset Q$ be the subset of such quadrangles with the additional condition that
\beq
|x+y| = n
\eeq
It will be convenient to use  the symbolic notation `$\squarediagup$'  for an element $(x, y, v, w)$ of $Q^-_n$;  and `$\ultriangle$' and `$\drtriangle$' for the corresponding triangles  $(x, y, z\!:=\!x+y)$ and $(z\!:=\! x+y, v, w)$. We also write `$\squarenodiag$' for a generic quadrangle in $Q$. 
Denote by $(V_{\ultriangle}, \Phi^g_{\ultriangle})$ and $(V_{\,\!\drtriangle},\Phi^g_{\!\drtriangle})$ the two 1-intertwiners $l\otimes m \to n$ and $n\otimes q\to r$.  
The composite 1-intertwiner $l \otimes m \otimes q \to r$ amounts to a SO(4) vector bundle ($W_{\squarediagup}, \Psi^g_{\squarediagup})$ over $Q^-_n$, with fibers $W_{\squarediagup}$ and maps $\Psi^g_{\squarediagup} \maps W_{\squarediagup} \to W_{g \squarediagup}$ given by the tensor products:
\beq \label{compo1int}
W_{\squarediagup}=V_{\ultriangle} \otimes  V_{\,\!\drtriangle}, \qquad \Psi^g_{\squarediagup} = \Phi^g_{\ultriangle} \otimes \Phi^g_{\!\drtriangle}
\eeq
If the two 1-intertwiners are irreducible, this is a line bundle: $W_{\squarediagup} \simeq \C$; and each of the maps $\Psi^g_{\squarediagup}$ acts by multiplication by a complex number. 

The composition of 1-intertwiners tells us also how to obtain an $\SO(4)$-invariant measure $\mu_{Q^-_n}$ on $Q^-_n$ from the measures $\mu_T$ and $\nu_{T'}$ on the two sets of triangles $(x, y, z)$ and $(z, v, w)$ given by (\ref{MeasInt}). It is defined by the formula: 
\beq \label{compo-meas}
\mu_{Q^-_n} = \pi \!\!\int \extd^3_nz \, \mu^z_{T} \otimes \nu^z_{T'}
\eeq
where $\extd^3_nz$ is the spherical measure (\ref{MeasSp}) and $\mu^z_{T}$ denotes the disintegration of $\mu_T$ with respect to $\extd^3_nz$, that is, $\extd\mu_T = \extd^3_nz \, \extd \mu_T^z$. The factor $\pi$ was inserted to simplify the formulas.  
Explicitly, $\extd \mu_{Q^-_n}$ reads: 
\beq \label{measQn}
\extd_l^3x \, \extd_m^3 y\, \extd_q^3 v \, \extd_r^3 w \, \delta(|x+y|^2-n^2) \, \delta^4(x+y+v-w) 
\eeq
The data (\ref{compo1int})  gives  a representation of SO(4) on the space
\beq \label{rep_n}
\direct \extd\mu_{Q^-_n}(\squarediagup)\, V_{\ultriangle} \otimes  V_{\,\!\drtriangle}
\eeq
of $L^2$ sections of the bundle $W_{\squarediagup}$.  
If the two 1-intertwiners are irreducible, $W_{\squarediagup}$ is line bundle; in this case the direct integral is just the space $L^2(Q^-_n, \mu_{Q^-_n})$  of square integrable functions on $Q^-_n$. 

Letting the label $n \!\in\! \R_+$ run within the range of values allowed by triangular inequalities, the union of the bundles over all subsets $Q^-_n \subset Q$ gives an SO(4) bundle over $Q$. We may consider the resulting SO(4) representation on $L^2$ sections with respect to the measure:
\beq
\mu_Q = \int \extd n^2 \mu_{Q^-_n}
\eeq
where $\extd n^2=2n \extd n$ is the Lebesgue measure on $\R_+$. Explicitly:
\beq \label{measQ}
\extd \mu_Q =  \extd_l^3x \, \extd_m^3 y\, \extd_q^3 v \, \extd_r^3 w \,  \delta^4(x+y+v-w) 
\eeq
This representation is the direct integral of the representations (\ref{rep_n}) labeled by $n$. For  $\squarenodiag \in Q$ and $\varphi \in  V_{\ultriangle} \otimes  V_{\,\!\drtriangle}$, the action of $g\in \SO(4)$ can be written as:
\beq \label{repQ1}
g (\squarenodiag, \varphi)  = (g\squarenodiag, (\Phi^g_{\ultriangle} \otimes  \Phi^g_{\!\drtriangle})(\varphi))
\eeq

In just the  same way, two 1-intertwiners
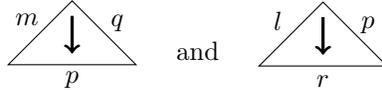
\[
\begin{array}{c}
\begin{tikzpicture}
\draw (0,0) -- (0.85,0.85) -- (1.7, 0) --cycle;
\path (0.25,0.6) node {\small $m$};
\path (1.45,0.55) node {\small $q$};
\path (0.85,-0.2) node {\small $p$};
\draw[->,very thick]
(0.85, 0.65) -- (0.85, 0.15);
\path (0.70, 0.37); 
\end{tikzpicture}
\end{array}
\quad \mbox{and} \quad 
\begin{array}{c}
\begin{tikzpicture}
\draw (0,0) -- (0.85,0.85) -- (1.7, 0) --cycle;
\path (0.25,0.6) node {\small $l$};
\path (1.45,0.55) node {\small $p$};
\path (0.85,-0.2) node {\small $r$};
\draw[->,very thick]
(0.85, 0.65) -- (0.85, 0.15);
\path (0.70, 0.37); 
\end{tikzpicture}
\end{array}
\]
can be composed as drawn:
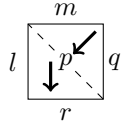
\[
\begin{array}{c}
\begin{tikzpicture}
\draw (0,0) rectangle (1, 1);
\path (-0.2, 0.5) node {\small $l$};
\path (0.5, 1.2) node {\small $m$};
\path (1.15, 0.5) node {\small $q$};
\path (0.5, -0.2) node {\small $r$};
\path (0.5, 0.5) node {\small $p$};
\draw[dashed] (1,0) -- (0.6,0.4);
\draw[dashed] (0.35, 0.65) -- (0, 1);
\path (0.53, 0.84);
\path (0.14, 0.37);
\draw[->, very thick]
(0.9,0.9) -- (0.6, 0.6);
\draw[->, very thick]
(0.3, 0.5) -- (0.3, 0.1);
\end{tikzpicture}
\end{array}
\]
Using an obvious extension of the above notations, we denote by $(V_{\dltriangle}, \Phi^g_{\dltriangle})$ and 
$(V_{\,\!\urtriangle}, \Phi^g_{\,\!\urtriangle})$ the two 1-intertwiners $l\otimes p \to r$ and $m\otimes q \to r$. The composite 1-intertwiner $l \otimes m \otimes q \to r$ amounts to a SO(4) vector bundle ($W_{\squarediagdown}, \Psi^g_{\squarediagdown})$ over the set $Q^+_p \subset Q$ of quadrangles (\ref{quadra}) such that
$|y+v| = p$,
where fibers $W_{\squarediagdown}$ and maps $\Psi^g_{\squarediagdown} \maps W_{\squarediagdown} \to W_{g \squarediagdown}$ are given by the tensor products:
\beq
W_{\squarediagdown} = V_{\dltriangle}\otimes  V_{\,\!\urtriangle}, \quad \Psi^g_{\squarediagdown} = \Phi^g_{\dltriangle}  \otimes \Phi^g_{\,\!\urtriangle}
\eeq
This data gives a representation of SO(4) on the space
\beq \label{rep_p}
\direct \extd\mu_{Q^+_p}(\squarediagup)\, V_{\dltriangle}\otimes  V_{\,\!\urtriangle}
\eeq
of $L^2$ sections on the bundle $W_{\squarediagdown}$,  where the measure $\mu_{Q^+_p}$ is defined by a formula analogous to (\ref{compo-meas}). Explicitly, $\extd\mu_{Q^+_p}$ reads:
\beq \label{measQn}
\extd_l^3x \, \extd_m^3 y\, \extd_q^3 v \, \extd_r^3 w \, \delta(|y+v|^2-p^2)  \, \delta^4(x+y+v-w) 
\eeq
Letting the label $p\in \R_+$  run within the range of values allowed by triangular inequalities, 
the direct integral of all representations (\ref{rep_p}) gives a SO(4) representation on $L^2$ sections of a bundle over $Q$, with respect to the measure:
\beq
\mu_Q = \int \extd p^2 \mu_{Q^+_p}
\eeq
also given by (\ref{measQ}). For  $\squarenodiag \in Q$ and $\psi \in  V_{\dltriangle}\otimes  V_{\,\!\urtriangle}$, the action of $g\in \SO(4)$ can be written as:
\beq \label{repQ2}
g (\squarenodiag, \psi)  = (g\squarenodiag, (\Phi^g_{\dltriangle}  \otimes \Phi^g_{\,\!\urtriangle})(\psi))
\eeq

The above constructions can be extended to more general triangulated surfaces with boundaries. 
For example,  by composing three 1-intertwiners as: 
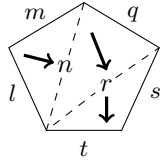
\[
\begin{array}{c}
\begin{tikzpicture}
\draw (-0.5, 0) --(-1, 1.1) -- (0, 1.7) -- (1,1.1) -- (0.5, 0) -- (-0.5, 0);
\draw[dashed] (0, 1.7) -- (-0.2,1);
\draw[dashed] (-0.3, 0.7) -- (-0.5, 0) -- (0.15, 0.5);
\draw[dashed] (0.42, 0.7) -- (1,1.1);
\path (-0.95, 0.5) node {\footnotesize $l$};
\path (-0.65, 1.55) node {\footnotesize $m$};
\path (0.65, 1.55) node {\footnotesize $q$};
\path (0.95, 0.5) node {\footnotesize $s$};
\path (0, -0.25) node {\footnotesize $t$};
\path (-0.25, 0.85) node {\small $n$};
\path (0.3, 0.6) node {\small $r$};
\draw[->, very thick]
(-0.8, 1) -- (-0.4, 0.9);
\draw[->, very thick]
(0.1, 1.3) -- (0.3, 0.8);
\draw[->, very thick]
(0.3, 0.45) -- (0.3, 0.1);
\end{tikzpicture}
\end{array}
\]
and by letting the labels $n,r\in \R_+$ on the dashed edges run within the range of values allowed by triangular inequalities, we obtain a 
representation of SO(4) on $L^2$ sections 
\beq \label{penta}
\direct \extd\mu_{\P}(\pentagon) \, V_{\ultriangle} \otimes  V_{\,\!\drtriangle} \otimes V_{\downtriangle} 
\eeq 
of a vector bundle over a set $\cP$ of pentagons $(x, y, v, w, z)$ of lengths $l,m,q, s,t$ in $\R^4$, endowed with the measure:
\beq \label{measP}
\extd\mu_\cP = \extd_l^3x \, \extd_m^3 y\, \extd_q^3 v \, \extd_s^3 w \, \extd_t^3 z\, \delta^4(x+y+v+w-z) 
\eeq
This representation is the direct integral of representations on $L^2$ sections of bundles over sets  $\cP_{nr}$ of pentagons with two fixed diagonal lengths, with respect to the measures $\mu_{\cP_{nr}}$ showing up in the decomposition:
\beq
\mu_\cP  = \int \extd n^2 \extd r^2 \, \mu_{\cP_{nr}}
\eeq
When the three 1-intertwiners are irreducible, the direct integral (\ref{penta}) is just the space $L^2(\cP, \mu_\P)$ of square integrable functions on $\cP$.  

\ss{2-Intertwiners on tetrahedra} \label{2int}

Given two composite 1-intertwiners
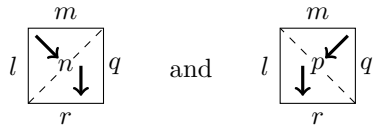
\[
\begin{array}{c}
\begin{tikzpicture}
\draw (0,0) rectangle (1, 1);
\path (-0.2, 0.5) node {\small $l$};
\path (0.5, 1.2) node {\small $m$};
\path (1.15, 0.5) node {\small $q$};
\path (0.5, -0.2) node {\small $r$};
\path (0.5, 0.5) node {\small $n$};
\draw[dashed] (0,0) -- (0.45,0.45);
\draw[dashed] (0.6, 0.6) -- (1, 1);
\path (0.45, 0.84);
\path (0.87, 0.37);
\draw[->, very thick]
(0.1,0.9) -- (0.4, 0.6);
\draw[->, very thick]
(0.7, 0.5) -- (0.7, 0.1);
\end{tikzpicture}
\end{array}\quad \mbox{and} \quad  \begin{array}{c}
\begin{tikzpicture}
\draw (0,0) rectangle (1, 1);
\path (-0.2, 0.5) node {\small $l$};
\path (0.5, 1.2) node {\small $m$};
\path (1.15, 0.5) node {\small $q$};
\path (0.5, -0.2) node {\small $r$};
\path (0.5, 0.5) node {\small $p$};
\draw[dashed] (1,0) -- (0.6,0.4);
\draw[dashed] (0.35, 0.65) -- (0, 1);
\path (0.53, 0.84);
\path (0.14, 0.37);
\draw[->, very thick]
(0.9,0.9) -- (0.6, 0.6);
\draw[->, very thick]
(0.3, 0.5) -- (0.3, 0.1);
\end{tikzpicture}
\end{array}
\]
constructed as above, a 2-intertwiner between these is drawn as:
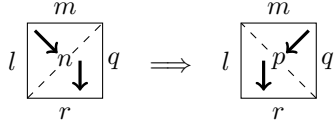
\[
\begin{array}{c}
\begin{tikzpicture}
\draw (0,0) rectangle (1, 1);
\path (-0.2, 0.5) node {\small $l$};
\path (0.5, 1.2) node {\small $m$};
\path (1.15, 0.5) node {\small $q$};
\path (0.5, -0.2) node {\small $r$};
\path (0.5, 0.5) node {\small $n$};
\draw[dashed] (0,0) -- (0.45,0.45);
\draw[dashed] (0.6, 0.6) -- (1, 1);
\path (0.45, 0.84);
\path (0.87, 0.37);
\draw[->, very thick]
(0.1,0.9) -- (0.4, 0.6);
\draw[->, very thick]
(0.7, 0.5) -- (0.7, 0.1);
\end{tikzpicture}
\end{array} \Longrightarrow  \begin{array}{c}
\begin{tikzpicture}
\draw (0,0) rectangle (1, 1);
\path (-0.2, 0.5) node {\small $l$};
\path (0.5, 1.2) node {\small $m$};
\path (1.15, 0.5) node {\small $q$};
\path (0.5, -0.2) node {\small $r$};
\path (0.5, 0.5) node {\small $p$};
\draw[dashed] (1,0) -- (0.6,0.4);
\draw[dashed] (0.35, 0.65) -- (0, 1);
\path (0.53, 0.84);
\path (0.14, 0.37);
\draw[->, very thick]
(0.9,0.9) -- (0.6, 0.6);
\draw[->, very thick]
(0.3, 0.5) -- (0.3, 0.1);
\end{tikzpicture}
\end{array}
\]
The two sides of this diagram correspond to the splitting of the boundary of a  tetrahedron with edge lengths $l,m,n,p,q,r$ into two pairs of triangles sharing an edge. 

Using the notations of the previous section, let $(V_{\ultriangle}, \Phi^g_{\ultriangle}), (V_{\,\!\drtriangle},\Phi^g_{\!\drtriangle}), (V_{\dltriangle}, \Phi^g_{\dltriangle})$ and  $(V_{\,\!\urtriangle}, \Phi^g_{\,\!\urtriangle})$ be the 1-intertwiners that label the four triangles. The two composite intertwiners $l\otimes m \otimes q \to r$
give two SO(4) vector  bundles $(V_{\ultriangle} \otimes  V_{\,\!\drtriangle}, \Phi^g_{\ultriangle} \otimes \Phi^g_{\!\drtriangle})$
and $(V_{\dltriangle}\otimes  V_{\,\!\urtriangle}, \Phi^g_{\dltriangle}  \otimes \Phi^g_{\,\!\urtriangle})$  over the sets $Q^{-}_n $ and $Q^+_p$, respectively. 
We also introduce the set 
\[Q_{np}= Q^-_n\cap Q^+_p\] of quadrangles $(x, y, v, w)$ with fixed diagonal lengths $|x+y|=n$ and $|y+v|=p$.  
Note that a generic element $\squarediag \in Q_{np}$ defines a tetrahedron of lengths $l, m, n, p, q,r$ embedded in $\R^4$. The natural action of SO(4) on $\R^4$ induces a transitive action on $Q_{np}$, which we write $(g, \squarediag)\mapsto g  \squarediag$. 
A 2-interwiner amounts to a family of linear maps
\beq \label{maps2int}
m^{}_{\squarediag} \maps V_{\ultriangle} \otimes  V_{\,\!\drtriangle} \to  V_{\dltriangle}\otimes  V_{\,\!\urtriangle}
\eeq
indexed by elements $\squarediag \in Q_{np}$, satisfying the intertwining rule:
\beq \label{2Int}
(\Phi^g_{\dltriangle}\otimes \Phi^g_{\urtriangle})  \, m^{}_{\squarediag } = 
 m^{}_{g \squarediag} \, (\Phi^g_{\ultriangle} \otimes  \Phi^g_{\!\drtriangle})  
\eeq
for all $g\in \SO(4)$.  

Since SO(4) acts transitively on $Q_{np}$,  this rule determines all maps $m^{}_{\squarediag}$ in terms of the map at a given point of $Q_{np}$. 
To fix the normalization of our 2-intertwiners, it  will be convenient to single out a `reference' quadrangle in  $Q_{np}$.  
For each set of values for the edge labels, we denote by $\squarediag^o$ be the unique quadrangle $(x_0, y_0, v_0, w_0)$ in $Q_{np}$ such that:
\beqa \label{reftet}
x_0 \!&=&\! l e_1, \nonumber \\ 
z_0 \!&\!:=&\! x_0+y_0 = n(\sin\gamma e_1 + \cos \gamma e_2),\nonumber \\
w_0 \!&=&\! r\cos\gamma' e_1 + r \sin \gamma'(\cos \theta e_2 + \sin\theta e_3)
\eeqa
for some angles $\gamma, \gamma', \theta \in [0,\pi]$, where $e_i$ are the basis vectors.  Thus, the triple of vectors $(x_0, z_0, w_0)$ defines a tetrahedron with positive orientation  in the  3-dimensional space ($e_1, e_2, e_3$). 
We also introduce the four triangles $(x_0, y_0, x_0+y_0)$, $(x_0+y_0, v_0, w_0)$, $(y_0, v_0, y_0+v_0)$ and $(x_0, y_0+v_0, w_0)$ induced by the reference tetrahedron; we use the symbols $\ultriangle(\squarediag^o)$,  $\drtriangle(\squarediag^o)$, $\urtriangle(\squarediag^o)$ and $\dltriangle(\squarediag^o)$ for these. Each of these triangles is the image by a unique SO(3) rotation
in the space $(e_1, e_2, e_3)$  of one of the reference triangles $\ultriangle^{\tiny{o}}$, $\drtriangle^{\tiny{o}}$, $\urtriangle^{\tiny{o}}$, and $\dltriangle^{\!\!\tiny{o}}$ lying in the plane $(e_1, e_2)$ and  specified by (\ref{reftriangle}).
We denote as $k_{\tinyultriangle}, k_{\tinydrtriangle}$, $k_{\tinyurtriangle}$ and $k_{\tinydltriangle}$ such SO(3) rotations:
\beqa \label{k's}
k_{\tinyultriangle}  \ultriangle^{\tiny{o}} &=& \ultriangle(\squarediag^o), \quad k_{\tinydrtriangle}\, \drtriangle^{\tiny{o}} = \drtriangle(\squarediag^o), \nonumber \\  
k_{\tinyurtriangle}\urtriangle^{\tiny{o}} &=& \urtriangle(\squarediag^o),\quad  k_{\tinydltriangle} \dltriangle^{\!\!\tiny{o}} =  \dltriangle(\squarediag^o) 
\eeqa
As a consequence of the rules (\ref{2Int}) and  (\ref{cocycle}), the map $m_{\squarediag}$  at any point $\squarediag\!=\!g\squarediag^o$  of  $Q_{np}$ can always be expressed as:
\beq \label{m0}
m^{}_{g \squarediag^o} = (\Phi^{g k_{\tinydltriangle}}_{\dltriangle^{\!\!\tiny{o}}} \otimes  \Phi^{g k_{\tinyurtriangle}}_{\urtriangle^{\tiny{o}}}) \,  m_0 \, \left(\Phi^{g k_{\tinyultriangle}}_{\ultriangle^{\tiny{o}}} \otimes \Phi^{g k_{\tinydrtriangle}}_{\drtriangle^{\tiny{o}}}\right)^{- 1}
\eeq
for some linear map $m_0 \maps V_{\ultriangle^{\!\tiny{o}}} \otimes V_{\drtriangle^{\tiny{o}}} \to V_{ \dltriangle^{\!\!\tiny{o}}} \otimes V_{\urtriangle^{\tiny{o}}}$ depending on the labels.   
Note that the map $m_0$  is  independent of our choice (\ref{reftet}) of positively oriented tetrahedron in the space  $(e_1, e_2, e_3)$. 
In fact, upon a rotation $\squarediag^o \to u \squarediag^o$ of the reference tetrahedron by $u\in \SO(3)$, the SO(3) rotations $k_{\!\vartriangle}$ defined in (\ref{k's}) simply change as 
$k_{\!\vartriangle} \to k^u_{\!\vartriangle} = u k_{\!\vartriangle}$. Hence such a rotation of the reference tetrahedron simply corresponds to a shift $g \to gu$ in $(\ref{m0})$.

When the 1-intertwiners  on the four triangles are all irreducible,  labeled by $s_i \in \Z$, $i=1,..4$,  all bundles are line bundles and the maps (\ref{maps2int}) act by multiplication by a complex number  $m^{}_{\squarediag} \in \C$.  In this case, the map $m_0$ showing up in  (\ref{m0}) only contributes to a normalization factor, which we set to one: 
\beq \label{norm}
 m_0 = 1
\eeq
Given such choices of a reference tetrahedron and normalization condition, the set of labels on the edges and triangles of a tetrahedron uniquely specifies a 2-intertwiner (\ref{maps2int}) via the formula (\ref{m0}). 
 
Note that our normalized 2-intertwiner depends on the orientation (here chosen to be positive) of the reference tetrahedron in the space $(e_1, e_2, e_3)$. In fact, a flip of the orientation by means of the U(1) rotation  $\squarediag^o \to h_{\pi} \squarediag^o$ of  angle $\pi$ around the plane $(e_1, e_2)$,  induces the change $k_{\!\vartriangle} \to h_{\pi} k_{\!\vartriangle} h_{\pi}$ of the rotations  (\ref{k's}).  Using the rule (\ref{rule}), this in turn leads to the rescaling $m_{\squarediag} \to (-1)^{\sum_i s_i} m_{\squarediag}$ of the 2-intertwiner. 

We also introduce the {\it dual}  2-intertwiner between the two composite 1-intertwiners 
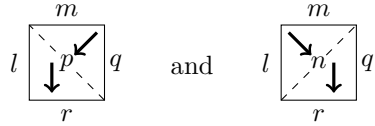
\[
\begin{array}{c}
\begin{tikzpicture}
\draw (0,0) rectangle (1, 1);
\path (-0.2, 0.5) node {\small $l$};
\path (0.5, 1.2) node {\small $m$};
\path (1.15, 0.5) node {\small $q$};
\path (0.5, -0.2) node {\small $r$};
\path (0.5, 0.5) node {\small $p$};
\draw[dashed] (1,0) -- (0.6,0.4);
\draw[dashed] (0.35, 0.65) -- (0, 1);
\path (0.53, 0.84);
\path (0.14, 0.37);
\draw[->, very thick]
(0.9,0.9) -- (0.6, 0.6);
\draw[->, very thick]
(0.3, 0.5) -- (0.3, 0.1);
\end{tikzpicture}
\end{array}\quad \mbox{and} \quad  \begin{array}{c}
\begin{tikzpicture}
\draw (0,0) rectangle (1, 1);
\path (-0.2, 0.5) node {\small $l$};
\path (0.5, 1.2) node {\small $m$};
\path (1.15, 0.5) node {\small $q$};
\path (0.5, -0.2) node {\small $r$};
\path (0.5, 0.5) node {\small $n$};
\draw[dashed] (0,0) -- (0.45,0.45);
\draw[dashed] (0.6, 0.6) -- (1, 1);
\path (0.45, 0.84);
\path (0.87, 0.37);
\draw[->, very thick]
(0.1,0.9) -- (0.4, 0.6);
\draw[->, very thick]
(0.7, 0.5) -- (0.7, 0.1);
\end{tikzpicture}
\end{array}
\]
drawn by reverting the double arrow as: 
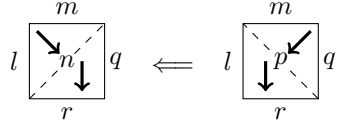
\[
\begin{array}{c}
\begin{tikzpicture}
\draw (0,0) rectangle (1, 1);
\path (-0.2, 0.5) node {\small $l$};
\path (0.5, 1.2) node {\small $m$};
\path (1.15, 0.5) node {\small $q$};
\path (0.5, -0.2) node {\small $r$};
\path (0.5, 0.5) node {\small $n$};
\draw[dashed] (0,0) -- (0.45,0.45);
\draw[dashed] (0.6, 0.6) -- (1, 1);
\path (0.45, 0.84);
\path (0.87, 0.37);
\draw[->, very thick]
(0.1,0.9) -- (0.4, 0.6);
\draw[->, very thick]
(0.7, 0.5) -- (0.7, 0.1);
\end{tikzpicture}
\end{array} \Longleftarrow  \begin{array}{c}
\begin{tikzpicture}
\draw (0,0) rectangle (1, 1);
\path (-0.2, 0.5) node {\small $l$};
\path (0.5, 1.2) node {\small $m$};
\path (1.15, 0.5) node {\small $q$};
\path (0.5, -0.2) node {\small $r$};
\path (0.5, 0.5) node {\small $p$};
\draw[dashed] (1,0) -- (0.6,0.4);
\draw[dashed] (0.35, 0.65) -- (0, 1);
\path (0.53, 0.84);
\path (0.14, 0.37);
\draw[->, very thick]
(0.9,0.9) -- (0.6, 0.6);
\draw[->, very thick]
(0.3, 0.5) -- (0.3, 0.1);
\end{tikzpicture}
\end{array}
\]
This is the unique 2-intertwiner whose maps 
\beq
\bar{m}_{\squarediag} \maps V_{\dltriangle}\otimes  V_{\,\!\urtriangle} \to  V_{\ultriangle} \otimes  V_{\,\!\drtriangle} 
\eeq
satisfy the normalization condition (\ref{norm}) for the {\sl flipped} reference tetrahedron $\bar{\squarediag}^o = h_{\pi} \squarediag^o$. 
The normalization is chosen in such a way that the maps $m_{\squarediag}$ of the 2-intertwiner and those $\bar{m}_{\squarediag}$ of its dual satisfy the identity: 
\beq \label{normalization}
(-1)^{\sum_i \!s_i} \, \mbox{Tr}[ \bar{m}_{\squarediag} \, m_{\squarediag} ]= 1
\eeq
where $\mbox{Tr}$ is the trace in $V_{\dltriangle}\otimes  V_{\,\!\urtriangle}$ and the sum is over the four spins $s_i\in \Z$ labeling the irreducible 1-intertwiners.

Just as 1-intertwiners define SO(4) representations, 2-intertwiners define SO(4) intertwining operators. 
To see this, let the labels $n, p \in \R_+$ on the dashed edges of the above diagrams run in the range of values allowed by the triangular inequalities.  The union of the bundles over all subsets $Q^-_n$ and $Q^+_p \subset Q$ give two SO(4) bundles over $Q$. Consider the map 
\beq \label{SO(4)int}
(\squarenodiag, \varphi)  \, \mapsto \,  (\squarenodiag, m^{}_{\squarediag}(\varphi))
\eeq
from sections of one bundle to sections of the other. In this formula, if `$\squarenodiag$' is the quadrangle $(x, y, v, w)$, `\squarediag'  symbolizes the corresponding element in the subset $Q_{np}$ with $n\!=\!|x+y|$ and $p\!=\!|y+v|$. The rule (\ref{2Int}) says that (\ref{SO(4)int}) defines an intertwiner between the two $\SO(4)$ representations (\ref{repQ1}) and (\ref{repQ2}).
Restricting to $L^2$ sections with respect to the measure (\ref{measQ}), we obtain a map
\beq \label{m}
M \maps \direct \extd\mu_Q(\squarenodiag) \, V_{\ultriangle} \otimes  V_{\,\!\drtriangle}\, \to \, \direct \extd\mu_Q(\squarenodiag)  V_{\dltriangle}\otimes  V_{\,\!\urtriangle}
\eeq 
This map is the direct integral of the maps (\ref{maps2int}).  
When the 1-intertwiners  on the four triangles are all irreducibles, $M$
is the diagonal operator on the space $L^2(Q, \mu_Q)$ of square integrable functions on $Q$ acting by multiplication by the function $\squarenodiag \mapsto m^{}_{\squarediag}$.

These constructions  can be extended to 2-intertwiners between more general composite 1-intertwiners. 
For example, if we supplement the six representations and the four 1-intertwiners considered above with two irreducible representations $s , t\in \R_+$ and a 1-intertwiner $(V_{\downtriangle}, \Phi_{\downtriangle})$ between $r\otimes s$ and $t$, we may consider the 2-intertwiner drawn as:
\vspace{0.2cm}
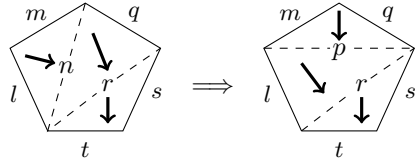
\[
\begin{array}{c}
\begin{tikzpicture}
\draw (-0.5, 0) --(-1, 1.1) -- (0, 1.7) -- (1,1.1) -- (0.5, 0) -- (-0.5, 0);
\draw[dashed] (0, 1.7) -- (-0.2,1);  \draw[dashed] (-0.3, 0.7) -- (-0.5, 0) -- (0.15, 0.5);
\draw[dashed] (0.42, 0.7) -- (1,1.1);
\path (-0.95, 0.5) node {\footnotesize $l$};
\path (-0.65, 1.55) node {\footnotesize $m$};
\path (0.65, 1.55) node {\footnotesize $q$};
\path (0.95, 0.5) node {\footnotesize $s$};
\path (0, -0.25) node {\footnotesize $t$};
\path (-0.25, 0.85) node {\small $n$};
\path (0.3, 0.6) node {\small $r$};
\draw[->, very thick]
(-0.8, 1) -- (-0.4, 0.9);
\draw[->, very thick]
(0.1, 1.3) -- (0.3, 0.8);
\draw[->, very thick]
(0.3, 0.45) -- (0.3, 0.1);
\end{tikzpicture}
\end{array}
\Longrightarrow 
\begin{array}{c}\begin{tikzpicture}
\draw (-0.5, 0) --(-1, 1.1) -- (0, 1.7) -- (1,1.1) -- (0.5, 0) -- (-0.5, 0);
\draw[dashed] (-1, 1.1) -- (-0.12, 1.1);
\draw[dashed] (0.2, 1.1) -- (1,1.1); 
\draw[dashed] (-0.5, 0) -- (0.15, 0.5);
\draw[dashed] (0.42, 0.7) -- (1,1.1);
\path (-0.95, 0.5) node {\footnotesize $l$};
\path (-0.65, 1.55) node {\footnotesize $m$};
\path (0.65, 1.55) node {\footnotesize $q$};
\path (0.95, 0.5) node {\footnotesize $s$};
\path (0, -0.25) node {\footnotesize $t$};
\path (0, 1.05) node {\small $p$};
\path (0.3, 0.6) node {\small $r$};
\draw[->, very thick]
(0, 1.6) -- (0, 1.2);
\draw[->, very thick]
(-0.5, 0.9) -- (-0.2, 0.5);
\draw[->, very thick]
(0.3, 0.45) -- (0.3, 0.1);
\end{tikzpicture}
\end{array}
\]
defined by taking  the tensor product of the maps (\ref{maps2int}) with identity maps $\unit_{\downtriangle} \maps V_{\downtriangle} \to V_{\downtriangle}$:
\beq \label{comp2int}
m^{}_{\squarediag} \otimes \unit_{\downtriangle}  \maps (V_{\ultriangle} \otimes  V_{\,\!\drtriangle}) \otimes V_{\downtriangle} \to  (V_{\dltriangle}\otimes  
V_{\,\!\urtriangle})\otimes V_{\downtriangle}
\eeq 
Letting the labels on the dashed lines edges run within the range
 of values allowed by the triangular inequalities, we obtain a map
\[ 
\direct \extd\mu_{\P}(\pentagon) \, V_{\ultriangle} \otimes  V_{\,\!\drtriangle}\otimes V_{\downtriangle}\, \to \, \direct \extd\mu_{\P}(\pentagon)  V_{\dltriangle}\otimes  V_{\,\!\urtriangle}\otimes V_{\downtriangle}
\]
that intertwines two SO(4) representations on $L^2$ sections of bundles over the set $\cP$ of pentagons of lengths $l,m,q,s,t$ in $\R^4$, endowed with the measure (\ref{measP}). We used our symbolic notation where $\pentagon$ denotes a generic pentagon $(x,y,v,w,z)$ in $\cP$; and $\squarediag$ denotes the corresponding quadrangle $(x,y,v, x+y+v)$ in $Q_{np}$ with $|x+y| = n$ and $|y+v| =p$. 
When the 1-intertwiners are all irreducibles, the vectors spaces are all one dimensional  and the maps (\ref{comp2int}) act  by multiplication by a complex number $m^{}_{\squarediag} \in \C$.  In this case, the SO(4) intertwiner is the diagonal operator on the space $L^2(\P, \mu_{\P})$ of square integrable functions on $\cP$ acting by multiplication by the function $\pentagon \mapsto m^{}_{\squarediag}$.

\s{10$j$ 2-symbols} \label{symbol}

We are now in a position to define and compute the weight that the state sum associates to each 4-simplex of the triangulation.
This weight, which we refer to as `10$j$ 2-symbol',  is a function of ten positive numbers labeling the irreducible representations of the Euclidean 2-group, and ten integer spins  labeling irreducible 1-intertwiners. As we have seen in the previous section, each set of such labels determines one normalized 2-intertwiner for each of the five  boundary tetrahedra of the 4-simplex. 
\ss{Definition}

10$j$ 2-symbols will be defined by taking an appropriate trace of the product of five normalized  2-intertwiners, depicted by the following diagram:
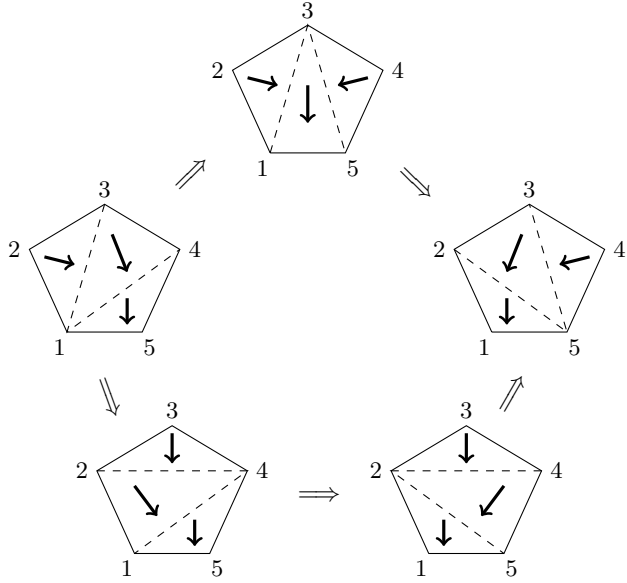
\[
\hspace{1cm}
\begin{array}{lllll}
\begin{array}{c}
\hspace{1.7cm}
\begin{tikzpicture}
\draw (-0.5, 0) --(-1, 1.1) -- (0, 1.7) -- (1,1.1) -- (0.5, 0) -- (-0.5, 0);
\draw[dashed] (-0.5, 0) -- (0, 1.7) -- (0.5, 0);
\path (-0.6, -0.2) node {\footnotesize $1$};
\path (-1.2, 1.1) node {\footnotesize $2$};
\path (0, 1.9) node {\footnotesize $3$};
\path (1.2, 1.1) node {\footnotesize $4$};
\path (0.6, -0.2) node {\footnotesize $5$};
\draw[->, very thick]
(-0.8, 1) -- (-0.4, 0.9);
\draw[->, very thick]
(0, 0.9) -- (0, 0.4);
\draw[->, very thick]
(0.8, 1) -- (0.4, 0.9);
\end{tikzpicture}
\end{array}
\vspace{-0.5cm} 
\\
\vspace{-0.4cm} \hspace{1.5cm} \rotatebox{45}{$\Longrightarrow$} \hspace{2.5cm}  \rotatebox{135}{$\Longleftarrow$}
\\
\begin{array}{cc}
\hspace{-1cm}
\begin{tikzpicture}
\draw (-0.5, 0) --(-1, 1.1) -- (0, 1.7) -- (1,1.1) -- (0.5, 0) -- (-0.5, 0);
\draw[dashed] (0, 1.7) -- (-0.5, 0) -- (1,1.1);
\path (-0.6, -0.2) node {\footnotesize $1$};
\path (-1.2, 1.1) node {\footnotesize $2$};
\path (0, 1.9) node {\footnotesize $3$};
\path (1.2, 1.1) node {\footnotesize $4$};
\path (0.6, -0.2) node {\footnotesize $5$};
\draw[->, very thick]
(-0.8, 1) -- (-0.4, 0.9);
\draw[->, very thick]
(0.1, 1.3) -- (0.3, 0.8);
\draw[->, very thick]
(0.3, 0.45) -- (0.3, 0.1);
\end{tikzpicture}
\hspace{2.5cm}
&
\begin{tikzpicture}
\draw (-0.5, 0) --(-1, 1.1) -- (0, 1.7) -- (1,1.1) -- (0.5, 0) -- (-0.5, 0);
\draw[dashed] (-1, 1.1) -- (0.5, 0) -- (0, 1.7);
\path (-0.6, -0.2) node {\footnotesize $1$};
\path (-1.2, 1.1) node {\footnotesize $2$};
\path (0, 1.9) node {\footnotesize $3$};
\path (1.2, 1.1) node {\footnotesize $4$};
\path (0.6, -0.2) node {\footnotesize $5$};
\draw[->, very thick]
(0.8, 1) -- (0.4, 0.9);
\draw[->, very thick]
(-0.1, 1.3) -- (-0.3, 0.8);
\draw[->, very thick]
(-0.3, 0.45) -- (-0.3, 0.1);
\end{tikzpicture}
\end{array}
\\
\vspace{-0.4cm} \hspace{0.5cm} \rotatebox{110}{$\Longleftarrow$} \hspace{5cm} {\rotatebox{60}{$\Longrightarrow$}}
\\
\hspace{-0.1cm}
\begin{array}{c}
\begin{tikzpicture}
\draw (-0.5, 0) --(-1, 1.1) -- (0, 1.7) -- (1,1.1) -- (0.5, 0) -- (-0.5, 0);
\draw[dashed] (-1, 1.1) --(1,1.1) -- (-0.5, 0);
\path (-0.6, -0.2) node {\footnotesize $1$};
\path (-1.2, 1.1) node {\footnotesize $2$};
\path (0, 1.9) node {\footnotesize $3$};
\path (1.2, 1.1) node {\footnotesize $4$};
\path (0.6, -0.2) node {\footnotesize $5$};
\draw[->, very thick]
(0, 1.6) -- (0, 1.2);
\draw[->, very thick]
(-0.5, 0.9) -- (-0.2, 0.5);
\draw[->, very thick]
(0.3, 0.45) -- (0.3, 0.1);
\end{tikzpicture}
\end{array}
\Longrightarrow 
\begin{array}{c}
\begin{tikzpicture}
\draw (-0.5, 0) --(-1, 1.1) -- (0, 1.7) -- (1,1.1) -- (0.5, 0) -- (-0.5, 0);
\draw[dashed]  (1,1.1) -- (-1, 1.1) --(0.5, 0);
\path (-0.6, -0.2) node {\footnotesize $1$};
\path (-1.2, 1.1) node {\footnotesize $2$};
\path (0, 1.9) node {\footnotesize $3$};
\path (1.2, 1.1) node {\footnotesize $4$};
\path (0.6, -0.2) node {\footnotesize $5$};
\draw[->, very thick]
(0, 1.6) -- (0, 1.2);
\draw[->, very thick]
(0.5, 0.9) -- (0.2, 0.5);
\draw[->, very thick]
(-0.3, 0.45) -- (-0.3, 0.1);
\end{tikzpicture}
\end{array}
\end{array}
\]

The details of the construction are as follows. Consider the standard 4-simplex $[12345]$, endowed with the standard orientation defined by the ordering  
$\{1,\cdots\!, 5\}$ of its vertices.  
The 4-simplex has boundary 
\beq \label{boundary}
[2345]-[1345]+[1245]-[1235]+[1234]
\eeq
where $[jklm]$ is the tetrahedron with vertices $\{j,k,l,m\}$. The sign indicates the induced orientation of the 4-simplex boundary: for the  tetrahedron $\hat \imath\!:=\! [jklm]$ opposite to the vertex $i$, the boundary orientation is defined as the ordering of its vertices in an even permutation of $(12345)$ where $i$ appears first. Each boundary tetrahedron shows up in the sequence (\ref{boundary}) with a (+) or a (-) sign depending on whether or not the numerical order agrees with the boundary orientation.

There is an irreducible representation for each edge $(ij)$, labeled by $l_{ij} \!\in\! \R_+$. The edges drawn as dashed lines on the diagram are coined `interior edges', those drawn as plain lines are coined `exterior edges'. For each triangle $[ijk]$ with vertices $i\!<\!j<\!k$,   there is an irreducible 1-intertwiner $l_{ij} \otimes l_{jk} \to l_{ik}$,  labeled by a spin $s_{ijk} \!\in\! \Z$. For each tetrahedron $\hat \imath\!=\! [jklm]$ with vertices $j\!<\!k<\!l\!<\!m$, there is a normalized 2-intertwiner $m^{\hat \imath}_{\squarediag}$ indexed by the vertex opposite to it in the 4-simplex;   we denote by  $\bar{m}^{\hat \imath}_{\squarediag}$ the corresponding dual 2-intertwiner.

Recalling the definitions of the previous sections,  $\cP$ is the set of pentagons $(x_{ij})$ of lengths $l_{ij}$ in $\R^4$, where $(ij)\in \{(12), (23), (34), (45), (15)\}$ runs over the five exterior edges.  
$\cP$ is endowed with the measure $\extd \mu_{\cP}$ given by the formula  (\ref{measP}): 
\beq \label{measP2}
\extd\mu_\cP = \prod_{{\mbox{\footnotesize ext.}} (ij)} \extd^3_{l_{ij}} x_{ij}  \, \delta^4(x_{12}+ .. +x_{45} - x_{15})
\eeq
where $\extd^3_{l_{ij}} x_{ij}$ is the spherical measure (\ref{MeasSp}) on the 3-sphere of radius $l_{ij}$; the delta function imposes the closure of the pentagons. 
Letting the labels on the interior edges run within the range of values allowed by the triangular inequalities, each of the five pentagonal figures of the diagram gives an SO(4) representation on $L^2$ sections of a vector bundle over $\cP$, defined as in (\ref{penta});  and each double arrow gives a SO(4) intertwiner between these, defined as in (\ref{comp2int}).  There is thus one such SO(4) intertwiner $M_{\hat \imath}$ for each tetrahedron $\hat \imath:= [jklm]$. We also denote by $\bar{M}_{\hat \imath}$ the SO(4) intertwiner determined by the dual 2-intertwiner.

We consider the product 
\beq
 M=\bar{M}_{\hat 2} \bar{M}_{\hat 4}  M_{\hat 1} M_{\hat 3}  M_{\hat 5}  
\eeq
of the five SO(4) intertwiners. The use of dual 2-intertwiners is dictated by the orientation of the tetrahedra: 
for each $\hat \imath$,  the product involves either $M_{\hat \imath}$ or the dual $\bar{M}_{\hat \imath}$, whether the  tetrahedron shows up in (\ref{boundary}) with a (+) or a (-) sign, that is,  whether $(ijklm)$ where $j\!<\!k\!<\!l\!<m$  is an even or an odd permutation of $(12345)$. 
By construction, $M$ is the direct integral over $\cP$ of maps $M_{\pentagon}$ indexed by elements $\pentagon \in \cP$ and given by:
\beq \label{Mp}
(\unit \otimes {\bar m}^{\hat 2}_{\squarediag_{\hat 2}}) ({\bar m}^{\hat 4}_{\squarediag_{\hat 4}}\otimes \unit)
 (\unit \otimes {m}^{\hat 1}_{\squarediag_{\hat 1}}) (\unit \otimes {m}^{\hat 3}_{\squarediag_{\hat 3}}) ({m}^{\hat 5}_{\squarediag_{\hat 5}}\otimes \unit) 
 \eeq
where $\squarediag_{\hat \imath}$ denotes the unique quadrangle (tetrahedron) in $\R^4$ that the point  $\pentagon \in \cP$ associates to $\hat \imath = [jklm]$.

10$j$ 2-symbols are defined by means of the identity:
\beq \label{2-10j}
\mbox{Tr} M =\kappa \!\int \prod_{{\mbox{\footnotesize int.}}(ij)}  \extd l_{ij}^2
\left
\{\begin{array}{ccc}
l_{e_1}&\cdots&l_{e_{10}} \\
s_{t_1}&\cdots&s_{t_{10}}
\end{array}
\right\} 
\eeq
where the integral is over the labels on the interior edges and $\kappa$ is an overall constant, which will later be chosen to be $\kappa =\pi^4/2^6$ for practical convenience. The 2-symbol, which we wrote here using brackets, depends on the labels $l_{e_j}$ on the ten edges $e_j$ and the labels  $s_{t_j}$ on the ten triangles $ t_j$ of the 4-simplex.  $\mbox{Tr} M$ denotes the {\it trace} of $M$, where the trace of a direct integral is defined as the integral of the trace:
\beq \label{trace}
\mbox{Tr} M := \int_{\cP} \extd\mu_\cP(\pentagon) \mbox{Tr} M_{\pentagon}
\eeq
The formula (\ref{2-10j}) should be understood as an identity of measures, as follows. Writing
$\cP = \sqcup_\ell \cP_\ell$ as the disjoint union of subsets $\cP_{\ell}$ of pentagons with given interior edge lengths $\mathbf{\ell}= (l_{ij})$, 
the identity (\ref{2-10j}) says that upon the decomposition 
\beq \label{muPl}
\mu_{\cP} = \kappa \! \int \prod_{{\mbox{\footnotesize int.}}(ij)} \extd l_{ij}^2 \, \mu_{\cP_{\mathbf{\ell}}}
\eeq
of the measure $\mu_{\cP}$ into measures $\mu_{\cP_{\mathbf{\ell}}}$ on $\cP_{\ell}$, the trace (\ref{trace})  decomposes into 10$j$ 2-symbols. 
Hence these symbols are explicitly given by the formula: 
\beq \label{def-symbols}
\left
\{\begin{array}{ccc}
l_{e_1}&\cdots&l_{e_{10}} \\
s_{t_1}&\cdots&s_{t_{10}}
\end{array}
\right\} 
:= \int_{\cP_\ell} \extd\mu_{\cP_\ell}(\pentagon) 
\begin{array}{c}
\begin{tikzpicture}[scale=0.7]
\draw
(0.3, 0.3) -- (1.2, 1.2);
\draw
(1.8, 1.2) -- (2.7, 0.3);
\draw
(2.7, -0.3) --(2.4, -1.2);
\draw
 (1.1, -1.5) -- (1.8, -1.5); 
\draw
(0.6, -1.2) -- (0.3 ,-0.3);
\draw
(0.7, -1.2) --(1.4, 1.1);
\draw
(2.3, -1.2) -- (1.6, 1.1);
\draw(0.4, 0) -- (0.9,0);
\draw (1.2, 0) --(1.8, 0);
\draw
(2.1, 0) --(2.6, 0);
\draw
(0.4, -0.15) -- (0.8, -0.4);
\draw (1.1, -0.6) -- (2, -1.2);
\draw (2.6, -0.15) --(2.2, -0.4);
\draw (1.9, -0.6) -- (1.6, -0.8);
\draw
(1, -1.2) -- (1.3, -1);
\path (0, 0) node {${m}^{\hat 5}_{\squarediag_{\hat 5}}$};
\path (1.5,1.5) node {${m}^{\hat 3}_{\squarediag_{\hat 3}}$};
\path (3.05, 0) node {${m}^{\hat 1}_{\squarediag_{\hat 1}}$};
\path (2.25, -1.5) node {$\bar  m^{\hat 4}_{\squarediag_{\hat 4}}$};
\path (0.75, -1.5) node {$\bar  m^{\hat 2}_{\squarediag_{\hat 2}}$};
\end{tikzpicture}
\end{array}
\eeq
where the diagram is the graphical representation of the trace $\mbox{Tr} M_{\pentagon}$ of the map (\ref{Mp}).  The integration  is over the subset $\cP_\ell \subset \cP$ of pentagons whose  edge lengths match with the labels  $l_{e_j}$. Note that  each of such pentagons defines a Euclidean 4-simplex embedded in $\R^4$.
 
\ss{Explicit computation}

The goal of the remainder of this section is to evaluate the integral in the right-hand side of (\ref{def-symbols}). 
This will enable us to write the 10$j$ 2-symbol as an explicit function of the labels. The result is the formula  (\ref{Form2-210j}) below.  

\sss{The measure} 

First, combining Equ. (\ref{measP2}), (\ref{muPl}) and (\ref{MeasSp}) gives the expression of the measure on $\cP_\ell$  in terms of the Lebesgue measure on $\R^4$. 
Solving the delta function in (\ref{measP2}) by integrating  over $x_{15}$ yields: 
\beq
\extd\mu_{\cP_\ell} = \frac{1}{\kappa \pi^5} \prod_{i=1}^4 \extd^4 x_{ii+1} \prod_{i<j} \delta(|x_{ij}|^2 - l^2_{ij}) 
\eeq
where we introduced the vectors $x_{ij}:=x_{ii+1} + ... + x_{j-1j}$. 
Each point of $\cP_\ell$ corresponds to a 4-simplex embedded in $\R^4$, with edge-vectors $x_{ij}$. 
Upon the action of SO(4), $\cP_\ell$ has two orbits, labeled by the orientation $\eta \!=\! \pm 1$ of the 4-simplices in $\R^4$.  The measure $\extd\mu_{\cP_\ell}$ being SO(4) invariant, it is thus equivalent to the product 
$\extd g \sum_{\eta=\pm 1}$ of a Haar measure on SO(4) and the counting measure on the set of orbits. The resulting Jacobian has been computed explicitly in \cite{BaratinFreidel3d}. For the value  $\kappa = \pi^9/2^6$, it leads to the identity:
\beq \label{measures}
\extd\mu_{\cP_\ell}  = \frac{1}{V(l_{ij})} \extd g \sum_{\eta=\pm 1}
\eeq
where $\extd g$ is the normalized Haar measure on SO(4) and $V(l_{ij})$ is $4!$ times the volume of a Euclidean 4-simplex with edge lengths $l_{ij}$.  

\sss{The trace}

The integrand in (\ref{def-symbols}) is the trace $\mbox{Tr} M_{\pentagon}$ of the map (\ref{Mp}).
By writing each of the five 2-intertwiners in the form (\ref{m0}), we obtain the formula: 
\beq \label{integrand}
\mbox{Tr} M_{\pentagon} = \langle \bigotimes_{\hat \imath} m^{\hat \imath}_{0} \,\,\, | \bigotimes_{\vartriangle= \langle  \hat\imath \hat\jmath  \rangle} 
\left(\Phi^{g_{\hat \jmath} k^{\hat \jmath}_{\vartriangle}}_{\vartriangle^{\!o}}\right)^{- 1} \Phi^{g_{\hat \imath} k^{\hat \imath}_{\vartriangle}}_{\vartriangle^{\!o}} \rangle
\eeq
whose right-hand side is the complex number obtained by tracing out the 1-intertwiner maps $\Phi^g_{\!\vartriangle^{\!o}}$ with the maps $m_0^{\hat \imath}$ showing up in the decomposition (\ref{m0}) of the 2-intertwiners. Our notation is as follows: 
given distinct $i,j,k,l,m$ with $k\!<\!l\!<m$, $\vartriangle=\langle  \hat\imath \hat\jmath  \rangle$ represents the triangle $[klm]$ common to the two tetrahedra $\hat \imath$ and $\hat \jmath$; moreover the pair $\langle  \hat\imath \hat\jmath  \rangle$ is ordered by requiring that $(ijklm)$ is an  even permutation of $(12345)$.  
The group elements $g_{\hat \imath}$ and $k^{\hat \imath}_{\vartriangle}$ are respectively the SO(4) and SO(3) rotations defined by:
\beq \label{defgk}
\squarediag_{\hat \imath}= g_{\hat \imath} \squarediag_{\hat \imath}^o, \qquad k^{\hat \imath}_{\!\vartriangle}\!\!\vartriangle^{\!\!o} \,\,=\,\, \vartriangle\!\!(\squarediag^o_{\hat \imath})
\eeq
when the tetrahedron $\hat \imath$  shows up in (\ref{boundary}) with a (+) sign ($\hat \imath = \hat 1, \hat 3, \hat 5$), and 
\beq \label{defgkdual}
\squarediag_{\hat \imath}= g_{\hat \imath} \bar{\squarediag}_{\hat \imath}^o, \qquad k^{\hat \imath}_{\!\vartriangle}\!\!\vartriangle^{\!\!o} \,\,=\,\, \vartriangle\!\!(\bar{\squarediag}^o_{\hat \imath})
\eeq
when $\hat \imath$ shows up in (\ref{boundary}) with a (-) sign ($\hat \imath = \hat 2, \hat 4$). We used the symbolic notation introduced in the previous section: $\squarediag_{\hat \imath}$ denotes the unique tetrahedron in $\R^4$ that the point  $\pentagon \in \cP$ associates to $\hat \imath = [jklm]$; $\squarediag^o_{\hat \imath}$ is the reference tetrahedron (\ref{reftet}) used to normalize 2-intertwiners; $\bar{\squarediag}^o_{\hat \imath} \!=\! h_\pi \squarediag^o_{\hat \imath}$, where $h_{\pi}$ is the rotation of angle $\pi$ around the plane $(e_1, e_2)$,  is the (flipped) reference tetrahedron used to normalize dual 2-intertwiners. Given a triangle $\vartriangle$ and a tetrahedron $\hat \imath$ adjacent to it, 
$\vartriangle^{\!\!o}$, $\vartriangle\!\!(\squarediag^o_{\hat \imath})$ (or $\vartriangle\!\!(\bar{\squarediag}^o_{\hat \imath})$) denote respectively the  reference triangle (\ref{reftriangle}) and the triangle that $\squarediag^o_{\hat \imath}$ (or $\bar{\squarediag}^o_{\hat \imath}$) associated to $\vartriangle$ in $\R^4$.
Note that,  written in the form (\ref{integrand}),  the trace of $M_{\pentagon}$  depends on the point $\pentagon \!\in\! \cP_\ell$  only through $g_{\hat \imath} \!:=\!g_{\hat \imath}(\pentagon)$.

Since in our case the 1-interwiners are irreducible, all vectors space are one-dimensional and all maps act by multiplication by a complex number. 
In this case, using the normalization (\ref{norm}) for the 2-intertwiners, the trace reduces to the following product: 
\beq \label{simpleTr} 
\mbox{Tr} M_{\pentagon} = \prod_{\vartriangle= \langle  \hat\imath \hat\jmath  \rangle} 
\left(\Phi^{g_{\hat \jmath} k^{\hat \jmath}_{\vartriangle}}_{\vartriangle^{\!o}}\right)^{- 1} \Phi^{g_{\hat \imath} k^{\hat \imath}_{\vartriangle}}_{\vartriangle^{\!o}}
\eeq
over the ten triangles $\vartriangle= \langle  \hat\imath \hat\jmath  \rangle$ of the 4-simplex. 

The next step is to observe that the contribution of the 1-intertwiner map corresponding to each triangle $\vartriangle= \langle  \hat\imath \hat\jmath  \rangle$ reduces to a phase $e^{is_{\!\vartriangle} \xi_{\!\vartriangle}}$, where $s_{\!\vartriangle} \in \mathbb{Z}$ is the spin label on  $\vartriangle$, and $\xi_{\!\vartriangle}$ is some angle in $[0,2\pi]$. 
Indeed, it is clear from the definitions (\ref{defgk}) and (\ref{defgkdual}) that the image of the reference triangle $\vartriangle^{\!\!o}$ by the rotation $g_{\hat \imath} k^{\hat \imath}_{\vartriangle}$ coincides with the triangle $\vartriangle\!\!\!(\squarediag_{\hat \imath})$ that $\squarediag_{\hat \imath}$ associates to $\vartriangle$ in $\R^4$. 
Since this triangle is common to $\hat \imath$ and $\hat \jmath$, we have that $\vartriangle\!\!\!(\squarediag_{\hat \imath}) = \vartriangle\!\!(\squarediag_{\hat \jmath})$. This shows that
the rotation  $(g_{\hat \jmath} k^{\hat \jmath}_{\vartriangle})^{-1} g_{\hat \imath} k^{\hat \imath}_{\vartriangle}$ stabilizes $\vartriangle^{\!\!o}$, hence belongs to U(1):
\beq \label{rotation}
(g_{\hat \jmath} k^{\hat \jmath}_{\vartriangle})^{-1} g_{\hat \imath} k^{\hat \imath}_{\vartriangle} \,\, := \,\, h_{\xi_{\!\vartriangle}}  \in \U(1)
\eeq          
for some $\xi_{\vartriangle} \in [0, 2\pi]$. Together with the relations (\ref{cocycle}) and (\ref{spincocycle}), it yields the formula:
\beq
\left(\Phi^{g_{\hat \jmath} k^{\hat \jmath}_{\vartriangle}}_{\vartriangle^{\!o}}\right)^{- 1} \Phi^{g_{\hat \imath} k^{\hat \imath}_{\vartriangle}}_{\vartriangle^{\!o}} 
 = 
\Phi^{(g_{\hat \jmath} k^{\hat \jmath}_{\vartriangle})^{-1} g_{\hat \imath} k^{\hat \imath}_{\vartriangle}}_{\vartriangle^{\!o}} = e^{is_{\!\vartriangle} \xi_{\!\vartriangle}}
\eeq

Note also that the rotations $h_{\xi_{\!\vartriangle}}\!\!:=\!h_{\xi_{\!\vartriangle}}\!(\pentagon)$ are all invariant under the SO(4) action $(g, \pentagon) \to g \pentagon$ on $\cP_\ell$. This means that, upon integration over $\cP_\ell$ with the measure (\ref{measures}), the angles $\xi_{\!\vartriangle}$, and hence the  integrand (\ref{def-symbols}), depend on the point $\pentagon \!\in \! \cP_\ell$ only through the orientation $\eta \! =\! \pm 1$ of the corresponding 4-simplex in $\R^4$.  
Hence, performing the integral of (\ref{simpleTr}) with respect to the measure (\ref{measures}) gives the quantity:
\beq \label{Form1-2-10j}
\frac{1}{V(l_{ij})} \sum_{\eta=\pm 1}\prod_{\vartriangle = \langle  \hat \imath \hat \jmath \rangle} e^{is_{\!\vartriangle} \xi_{\!\vartriangle}\!(\eta)}
\eeq

\sss{Relation to dihedral angles} 

The last step is to relate the angles $\xi_{\!\vartriangle}$ to the dihedral angles of the 4-simplex. The (interior) 
dihedral angle $\phi_{\!\vartriangle}\in [0,\pi]$ between the two tetrahedra sharing the triangle $\vartriangle = \langle  \hat \imath \hat \jmath \rangle$  is defined as ($\pi$ minus) the angle between their outwards unit normal vectors $n_\imath$ and $n_{\jmath}$: 
\beq
\cos \phi_{\!\vartriangle} = -n_\imath \cdot n_\jmath
\eeq

It is clear that the two angles $\xi_{\!\vartriangle}$ and $\phi_{\!\vartriangle}$ are closely related. For example we may note that, since by construction $k^{\hat \imath}_{\vartriangle}$ leaves $e_4$ invariant and $g_{\hat \imath}$ maps it to a vector normal to the tetrahedron $\imath$, the image of $e_4$ by 
$g_{\hat \imath} k^{\hat \imath}_{\vartriangle}$ must coincide with $\pm n_\imath$. This means the scalar product $e_4 \cdot h_{\xi_{\!\vartriangle}}e_4$ equals $ n_\imath \cdot n_\jmath$ modulo a sign, and thus
\beq
|\cos\xi_{\!\vartriangle}| =  |\cos \phi_{\!\vartriangle}|
\eeq 
which says $\xi_{\!\vartriangle} = \pm \phi_{\!\vartriangle}$ or $\pi \pm \phi_{\!\vartriangle}$.

The exact relation is the following: let $\vartriangle = \langle  \hat \imath \hat \jmath \rangle$ represent the triangle $[klm]$ common to the tetrahedra $\hat \imath$ and $\hat \jmath $, where $k\!<\!l\!<m$ and $(ijklm)$ is an even permutation of $(12345)$. 
The angle $\xi_{\!\vartriangle}$ of the rotation (\ref{rotation}) is  given by
\beq \label{relation-angle}
\xi_{\!\vartriangle} = \pi + \eta \phi_{\!\vartriangle}
\eeq
where $\eta$ is the orientation of the 4-simplex in $\R^4$, that is,  if $(x_{ij})_{i<j}$ are the edge vectors, the sign of the determinant $\det(x_{12}, \cdots, x_{15})$.

There are various ways to prove the relation (\ref{relation-angle}). An elegant algebraic proof relies on the following lemma. Let $a_{ij}^{kl}$ be the rotation of angle $\theta_{ij}^{kl}$ around the plane $(e_1, e_4)$, where  $\theta_{ij}^{kl}\in [0,\pi]$ is the 3-dimensional dihedral angle between the faces $[ijk]$ and $[ijl]$ in the tetrahedron $\hat m$.  
%
\begin{lemma} 
Given any permutation $(ijklm)$ of $(12345)$, the following equation for the triple of angles $(\xi_{ijk},  \xi_{ijl}, \xi_{ijm})$ in $[0,2\pi]^3$:
\beq \label{i}
h_{\xi_{ijk}} a_{ij}^{kl} \, h_{\xi_{ijl}} a_{ij}^{lm} \, h_{\xi_{ijm}}  a_{ij}^{mk} = \unit
\eeq
has exactly two solutions $\xi^\pm_{ijk}$ given by:
\be \label{sol-angles}
\xi^\pm_{ijk} = \pi \pm  \phi_{ijk}
\ee
where $\phi_{ijk}$ is the dihedral angle between the two tetrahedra sharing the triangle $[ijk]$.
\end{lemma}
%
%
%
There is one such identity (\ref{i}) satisfied by the dihedral angles (\ref{sol-angles}) for each edge $(ij)$ of the 4-simplex.  
These identities can be understood as vanishing curvature conditions around the edges of the 4-simplex \cite{BaratinFreidel3d, BiancaSimone}.

We refer to the Appendix B of \cite{BaratinFreidel3d}  for the proof of a {\it 3-dimensional} analogue of this Lemma. One dimension down, 
the analogues of the equations (\ref{i}) hold in SO(3) and are associated to the vertices of a tetrahedron $[iklm]$ in $\R^3$. 
In the usual basis $(e_2, e_3, e_4)$ of $\R^3$, they take the form:
\beq
 h_{\xi_{ik}}  a_{i}^{kl} h_{\xi_{il}} a_{i}^{lm} h_{\xi_{im}} a_{i}^{mk}  = \unit
\eeq
where $h_{\xi}$ is here the rotation of axis $e_2$ and angle $\xi$,  and $a_i^{kl}$ is the rotation of axis $e_4$ and angle $\theta_i^{kl}$, where $\theta_i^{kl} \in [0, \pi]$ is the angle between the edges $(ik)$ and $(il)$. 
These equations have two sets of solutions $\xi_{ik}^\pm \in [0,  2\pi]$ related to the 3-dimensional dihedral angles $\phi_{ik}$ as:
\beq \label{3dresult}
\xi_{ik}^\pm = \pi \pm \phi_{ik}
\eeq
We can summarize the correspondence between the 3d and 4d cases as follows: 
\begin{center}
{\small
\begin{tabular}{c|c}                   
4d            &   3d    
\rule[1.4em]{0em}{0em} \rule[-.8em]{0em}{0em}   \\      \hline
basis $(e_1, e_2, e_3, e_4)$  & basis $(e_2, e_3, e_4)$ \\ 
4-simplex $[ijklm]$    &  tetrahedron $[iklm]$  \\                
$a_\theta$ around $(e_1,e_4)$-plane  & $a_{\theta}$ around axis $e_4$ \\
$h_{\phi}$ around $(e_1, e_2)$-plane & $h_{\phi}$ around axis $e_2$ \\
3d dihedral  angles $\theta_{ij}^{kl}$   &  2d angles    $\theta_{i}^{kl}$             \\     
4d dihedral angles $\phi_{ijk}$  &  3d dihedral angles $\phi_{ik}$ \\
SO(4) edge identities  & SO(3) vertex identities 
\rule[-.8em]{0em}{0em}     \\     \hline
\end{tabular}} \vskip 1em
\end{center}
\vskip 0.5em

The proof of Lemma 1 can be inferred from its 3-dimensional analogue,  by considering the orthogonal projection $P \maps \R^4 \to \R^3$ onto the 3d space $(e_2, e_3, e_4)$ orthogonal to the basis vector $e_1$. Assuming the vertex $i$ of the 4-simplex $[ijklm]$ is at the origin and the edge $(ij)$ is along $e_1$, $P$ maps the 4-simplex to a tetrahedron $[iP(k)P(l)P(m)]$. The key observation is that the 3d dihedral angle $\theta_{ij}^{kl}$ between the faces $[ijk]$ and $[ijl]$ of the 4-simplex equals the 2d angle $\theta_i^{kl}$ between the edges $(iP(k))$ and $(iP(l))$ of the image tetrahedron; and the 4d dihedral angle $\phi_{ijk}$ equals  the 3d dihedral angle $\phi_{ik}$ between the faces $[iP(k)P(l)]$ and $[iP(k)P(m)]$ of the image tetrahedron.  
In short, the projection $P$ maps the left column of the above table to the right one. Using this correspondence, the Lemma 1 can be immediately deduced from the 3-dimensional result (\ref{3dresult}).

To see why the angles (\ref{rotation}) satisfy the relations (\ref{i}) of the Lemma, pick an edge $(ij)$ of the 4-simplex and let 
$(ijklm)$ be an even permutation of $(12345)$. The triangle common to the tetrahedra $\hat l$ and $\hat m $ 
is represented by the symbol $\vartriangle = \!\langle  \hat l \hat m \rangle$ (resp. \!\!$\vartriangle = \!\langle  \hat m \hat l \rangle$) if $i,j,k$ are numerically ordered as $\sigma(i)\!<\!\sigma(j) \!<\!\sigma(k)$ by an even  (resp.\!\! odd) permutation $\sigma$.  The reference triangle $\vartriangle^{\!\!o}$, whose edge vector $x_{\sigma(i)\sigma(j)}$ is along $e_1$ and whose edge vector $x_{\sigma(i)\sigma(k)}$ is in the plane  $(e_1, e_2)$ and points in the direction of $e_2$, is the image by some SO(3) rotation $\sigma_{ijk}$ of the triangle in $\R^4$ whose edge vector $x_{ij}$ is along $e_1$ and whose edge vector $x_{ik}$ is in the plane  $(e_1, e_2)$ and points in the direction of $e_2$. For example if $\!j\!<\!i\!<\!k$, i.e if $\sigma$ simply swaps $i$ and $j$, the action of $\sigma_{ijk}$ can be drawn as:
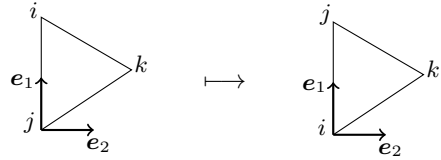
\[
\begin{array}{c}
\begin{tikzpicture}
\draw (0,0)--(0, 1.5)--(1.2 ,0.8)--(0,0);
\draw[->, thick] (0,0)-- (0,0.7);
\draw[->, thick] (0,0)--(0.7,0);
\path (-0.25,0.65) node {\footnotesize $\bm e_1$};
\path (0.76,-0.2) node  {\footnotesize $\bm  e_2$};
\path (-0.15, 0.1) node {\footnotesize $j$};
\path (-0.1,1.6) node {\footnotesize $i$};
\path  (1.33,0.88) node {\footnotesize $k$};
\end{tikzpicture}
\end{array}\quad  \longmapsto  \quad
\begin{array}{c}
\vspace{-0.1cm}
\begin{tikzpicture}
\draw (0,0)--(0, 1.5)--(1.2 ,0.8)--(0,0);
\draw[->, thick] (0,0)-- (0,0.7);
\draw[->, thick] (0,0)--(0.7,0);
\path (-0.25,0.65) node {\footnotesize $\bm e_1$};
\path (0.76,-0.2) node  {\footnotesize $\bm  e_2$};
\path (-0.15, 0.1) node {\footnotesize $i$};
\path (-0.1,1.6) node {\footnotesize $j$};
\path  (1.33,0.88) node {\footnotesize $k$};
\end{tikzpicture}
\end{array}
\]
Note that if $\sigma$ is  even, the two triangles have the same orientation in the plane $(e_1, e_2)$ and $\sigma_{ijk}$ leaves the whole plane $(e_3, e_4)$ invariant. If $\sigma$ is odd, the two triangles have opposite orientations and $\sigma_{ijk} \!=\! \tilde{\sigma}_{ijk} a_{\pi}$, where $\tilde{\sigma}_{ijk}$
leaves the whole plane $(e_3, e_4)$ invariant and $a_\pi$ is the rotation of angle $\pi$ around the plane $(e_1, e_4)$. 

If we then let $k_{ijk}^{\hat m}:=k_{\vartriangle}^{\hat m}$,  the rotations 
\beq \label{newh}
h_{\xi_{ijk}}= \sigma^{-1}_{ijk} (g_{\hat{l}} k_{ijk}^{\hat{l}})^{- 1} g_{\hat{m}} k_{ijk}^{\hat{m}} \sigma_{ijk}
\eeq
coincide with those defined in (\ref{rotation}). This is clear when $\sigma$ is even: in this case $\vartriangle = \!\langle  \hat l \hat m \rangle$ and  $\sigma_{ijk}$ commutes with all U(1) elements. When 
$\sigma$ is odd, then $\vartriangle = \!\langle  \hat m \hat l \rangle$ and $\sigma_{ijk} \!=\!\tilde{\sigma}_{ijk} a_{\pi}$ where $\tilde{\sigma}_{ijk}$ commutes with all U(1) elements, and we can conclude by using the equality
$a_{\pi}  h^{-1}_{\chi_{\vartriangle}} a_\pi = h_{\chi_{\vartriangle}}$. 
Observe also that the rotations
\beq \label{newa}
a_{ij}^{kl}=\sigma_{ijk}^{-1} (k_{ijk}^{\hat m})^{-1} k_{ijl}^{\hat m} \sigma_{ijl} 
\eeq
are those of the Lemma, of angle $\theta_{ij}^{kl}$ around the plane $(e_1, e_4)$. To see this,  let the edge vector $x_{ij}$ be along $e_1$, the edge vector $x_{il}$ be in the plane $(e_1, e_2)$ and point in the direction of $e_2$,  and the tetrahedron $\hat m$ be in the hyperplane $(e_1, e_2, e_3)$ with the orientation of its reference tetrahedron; since $(ijklm)$ is an even permutation of $(12345)$,  this means that $x_{ik}$ points in the direction opposite to $e_3$. The action of $a_{ij}^{kl}$ places the edge vector $x_{ik}$ in the plane $(e_1, e_2)$,  pointing in the direction of $e_2$:
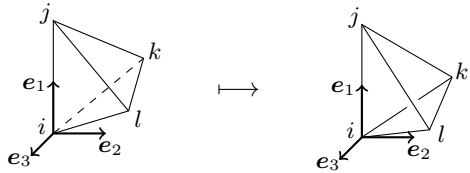
\[
\begin{array}{c}
\begin{tikzpicture}
\draw (0,0)-- (0,1.5) -- (1,0.3)--(0,0);
\draw (0,1.5)--(1.2 ,1)--(1,0.3);
\draw[dashed] (0,0)--(1.2 ,1);
\draw[->, thick] (0,0)-- (0,0.7);
\draw[->, thick] (0,0)--(0.7,0);
\draw[->, thick] (0,0)--(-0.3,-0.3);
\path (-0.25,0.65) node {\footnotesize $\bm e_1$};
\path (0.76,-0.2) node  {\footnotesize $\bm  e_2$};
\path (-0.45, -0.35) node {\footnotesize $\bm  e_3$};
\path (-0.15, 0.1) node {\footnotesize $i$};
\path (-0.1,1.6) node {\footnotesize $j$};
\path  (1.13,0.2) node {\footnotesize $l$};
\path (1.35 ,1.1)node {\footnotesize $k$};
\end{tikzpicture}
\end{array}\quad  \longmapsto  \quad
\begin{array}{c}
\vspace{-0.1cm}
\begin{tikzpicture}
\draw (0,0)--(0, 1.5)--(1.2 ,0.8);
\draw(0,0)--(0.6,0.4);
\draw (0.75, 0.5)--(1.2,0.8);
\draw (0,1.5)--(0.9,0.1)--(1.2 ,0.8);
\draw (0.9,0.1)--(0,0);
\draw[->, thick] (0,0)-- (0,0.7);
\draw[->, thick] (0,0)--(0.7,0);
\draw[->, thick] (0,0)--(-0.3,-0.3);
\path (-0.25,0.65) node {\footnotesize $\bm e_1$};
\path (0.76,-0.2) node  {\footnotesize $\bm  e_2$};
\path (-0.45, -0.35) node {\footnotesize $\bm  e_3$};
\path (-0.15, 0.1) node {\footnotesize $i$};
\path (-0.1,1.6) node {\footnotesize $j$};
\path  (1.33,0.88) node {\footnotesize $k$};
\path (1.06,0.05) node {\footnotesize $l$};
\end{tikzpicture}
\end{array}
\]

The rotations (\ref{newh}) are thus solutions of the equation (\ref{i}). 
Applying the lemma yields $\xi_{\!\vartriangle} \!=\! \pi + \epsilon \phi_{\!\vartriangle}$ for some $\epsilon \!=\! \pm1$ which does not depend on the triangle $\vartriangle$. It is then straightforward to relate  the sign $\epsilon$ to the orientation $\eta$ of the 4-simplex in $\R^4$ and reach the conclusion (\ref{relation-angle}).

Plugging the result (\ref{relation-angle}) into the formula (\ref{Form1-2-10j}) and summing over the orientation label gives the final expression for the 2-10$j$ symbol: 
\beq \label{Form2-210j}
\left
\{\begin{array}{ccc}
l_{e_1}&\cdots&l_{e_{10}} \\
s_{t_1}&\cdots&s_{t_{10}}
\end{array}
\right\}  = 
(-1)^{\sum_t \! s_t} \, \frac{\cos \left[\sum_{t}\! s_t \phi_t(l_e)\right]}{V(l_e)} 
\eeq
where both sums on the right-hand side are over the ten triangles of the 4-simplex; and $\phi_t$ is the dihedral angle between the two tetrahedra sharing the triangle $t$. 

\ss{State sum model}

Let $\Delta$ be a triangulated orientable 4-manifold. With any assignment of an irreducible representation $l_e\in \R_+$ of the Euclidean 2-group to each edge $e$ and an irreducible 1-intertwiner $s_t \in \Z$ to each triangle $t$ of $\Delta$, the model associates a weight $W_{\Delta}(l_e, s_t) \in \R$, given by the formula:
\beq \label{final-weight}
W_{\Delta} = \prod_{t } 2 A_t(l_e) \prod_{\tau } (-1)^{\sum_{i} \!s_{t_i}}  \,  \prod_{\sigma }\, 
\left\{\begin{array}{ccc}
l_{e_1}&\cdots&l_{e_{10}} \\
s_{t_1}&\cdots&s_{t_{10}}
\end{array}
\right\} 
\eeq
The products are over the triangles $t$, the tetrahedra $\tau$ and the 4-simplices $\sigma$ of the manifold.  
In the case of a manifold with boundary, the triangles and tetrahedra are those in the interior, i.e which do not lie on the boundary. 

The weight for each triangle corresponds to the volume of the measures (\ref{MeasInt}): it is equal to twice the area of the triangle  with edge lengths $l_e$ when the triangular inequality is satisfied, and zero otherwise. The weight for each tetrahedron corresponds to the normalization (\ref{normalization}) of the  2-intertwiners: in general it is defined to be  $1/\mbox{Tr}[ \bar{m}_{\squarediag} \, m_{\squarediag} ]$, which insures that $W_{\Delta}$ is independent of the normalization choice for an orientable manifold. For the normalization (\ref{norm}), it reduces to a sign factor $(-1)^{\sum_i \!s_{t_i}}$, where the sum is over the four triangles of the tetrahedron. These signs can be absorbed into the 4-simplex weight, giving a factor $i^{2 \sum_t \! s_t }$ which compensates the signs showing up in (\ref{Form2-210j}). Hence the formula agrees with (\ref{weights}) for a closed manifold; in general the two formulas agree up to a global sign depending only on the boundary data. 

The partition function is formally  defined by summing up these weights over all values of the labels, using the Lebesgue measures $\extd l_e^2$ for the real variables $l_e \!\in\! \R_+$ and the counting measure $\frac{1}{2\pi} \!\sum_{s_t}$ for the integer variables $s_t \!\in\! \Z$. Thus the 2-categorical state sum model constructed here is the same as the model of \cite{BaratinFreidel}. 

The physical interpretation of the model is best understood by writing the weight $W_{\Delta}$ in terms of the exponential of a discrete classical action, using the expression (\ref{Form1-2-10j}) of the 10$j$ 2-symbols. This action, which depends on the labels $l_e, s_t$ and an orientation $\eta_\sigma\!=\!\pm 1$  for each 4-simplex, reads
\beq \label{action} 
S_{\Delta}(l_e, s_t; \eta_\sigma) = \sum_t s_t \, \omega_t(l_e, \eta_\sigma)
\eeq  
where the sum is over all the triangles $t$. If we regard the edge labels $\{l_e\}$ as defining a discrete geometry in the sense of Regge calculus, 
the functions $\omega_t(l_e, \eta_\sigma)$, defined as  $\omega_t \! =\!\sum_{\sigma\supset t} \eta_\sigma \phi_t^\sigma$,  are the deficit angles associated to the triangles: they represent curvature in this geometry. The equations of motion obtained by varying the action (\ref{action}) imply that the deficit angles are trivial, i.e that the geometry is flat. 

The topological invariance of the state sum has been discussed in detail in \cite{BaratinFreidel}. The core of this discussion is an hexagonal identity  satisfied by the weight (\ref{final-weight}), giving a four-dimensional analogue of the Biedenhard-Elliot identity of 6$j$ symbols and an algebraic expression of the four-dimensional Pachner move invariance. Upon a regularization procedure resulting from the gauge-fixing of symmetries of the action (\ref{action}), it was shown in \cite{BaratinFreidel} that the partition function of the model reproduces the formula for the 4-manifold invariant defined by Korepanov in \cite{Korepanov1, Korepanov2, Korepanov3}.

\section{Outlook} \label{outlook}

We have developed in explicit detail a state sum model starting from the 2-category of representations of a 2-group. 
We have defined and computed the 4-simplex weights, which may be viewed as a categorified analogue of Racah-Wigner 6$j$ symbols. 
The 2-group we considered is a categorifed version of the four-dimensional Euclidean group ISO(4): though it is built from the same ingredients as the Euclidean group,  it differentiates the roles of rotations and translations by treating the former as objects in a category and the latter as morphisms. 
As anticipated in \cite{BarrettMackaay}, the resulting model has a remarkable geometrical flavour, where each set of irreducible representations labeling the edges of the triangulations are interpreted as specifying a Euclidean geometry on the underlying manifold. 

Our construction bridges results  from several works in the recent literature. First, it gives an explicit realization of the proposal of \cite{BarrettMackaay, CraneSheppeard} to use 2-group representations to define new state sum models in four dimensions.  Second, as shown in \cite{BaratinFreidel}, it provides a state sum formulation of the  4-manifold invariant constructed in \cite{Korepanov1, Korepanov2, Korepanov3}. Third, it unravels the algebraic structure underlying the state sum formulation of the Feynman amplitudes of quantum field theory on flat spacetime proposed in 
\cite{BaratinFreidel}. 

Our results suggest several interesting avenues for future work. To clarify the physical meaning of the model, an important step will be to identify the corresponding classical field theory (if any). A natural candidate is a higher gauge generalization of BF theory called `BFCG theory',  involving flat 2-connections \cite{BFCG}. We expect the model (\ref{final-weight}) to provide a state sum formulation of such a theory having the Euclidean 2-group as gauge 2-group. 

State sum models built from group representations have well-known generalizations obtained by replacing the group by a quantum group.  
In the case of the Lie group SU(2), one such generalization based on the quantum deformation $U_q(\mathrm{sl}_2)$ for $q$ root of unity  leads to the Turaev-Viro and Crane-Yetter models, which are finite and produce genuine manifold invariants \cite{TuraevViro, CKY}. It is worth investigating the construction of analogous models in the context of 2-groups, using  a suitable notion of quantum 2-group  \cite{Majid-2gp}. In the case of the Euclidean 2-group, this may lead to a regularization of the state sum (\ref{final-weight})  alternative  to the one provided by gauge fixing. Moreover, from the point of view  of \cite{BaratinFreidel}, this may enable one to propose and study possible dimension-full deformations of ordinary quantum field theory. This strategy has already revealed particularly useful to understand the coupling of matter fields to three-dimensional quantum gravity in the context of spin foam models  \cite{Barrett, PR3a, PR3b, BaratinFreidel3d}. 

Finally, several works pointed out the possible relevance of higher algebraic structures, and in particular the Poincar\'e 2-group, for the formulation of a model of quantum gravity in four dimensions \cite{BarrettMackaay, CraneSheppeard, Baratin-Wise, Mikovic, BaezWise}. 
Clarifying this relationship is yet another area for future study.

\vspace{0.7cm}

{\bf \large Acknowledgements}

 \vspace{0.2cm}
 
It is a pleasure to thank John Baez and Derek Wise for an extensive and fruitful collaboration on the representation theory of 2-groups. We also thank Jeffrey Morton and Hendryk Pfeiffer for insightful exchanges on this topic.  A.B gratefully acknowledges support of the Humboldt Foundation through a Feodor Lynen Fellowship. Research at Perimeter Institute is supported by the Government of Canada through Industry Canada and by the Province of Ontario through the Ministry of Research and Innovation. 

\end{document}